\documentstyle[12pt]{article}
\def \qed{\hfill $\vrule height 2.5mm  width 2.5mm depth 0mm $}

\newtheorem{th}{Theorem}
\newtheorem{pr}{Proposition}
\newtheorem{col}{Corollary}
\newtheorem{de}{Definition}

\def\al{\alpha}
\def\Res{\hbox{\rm Res}}
\def\mod{\hbox{\rm mod~}}
\def\Fl{\hbox{$Fl_n$}}
\def\ds{\displaystyle}
\def\wh{\widehat}
\def\wt{\widetilde}
\def\ld{\lambda}
\def\covec{\buildrel\leftarrow\over }
\font\germ=eufm10

\def\s{\hbox{\germ S}}

\def\a{\hbox{\germ A}}
\def\t{\hbox{\germ T}}
\begin{document}
\title{Quantum double Schubert polynomials, quantum Schubert polynomials
and Vafa--Intriligator formula.}
\vskip 0.4cm
\author{\Large {Anatol N. Kirillov}\thanks{On leave from {\it Steklov 
Mathematical Institute,
Fontanka 27, St.Petersburg, 191011, Russia}}~ \Large{and Toshiaki 
Maeno}\thanks{Supported by JSPS Research Fellowships for Young Scientists}
\medskip \\ 
{\small {\it Department of Mathematical Sciences, University of Tokyo,}}\\
{\small {\it Komaba, Meguro-ku, Tokyo 153, Japan}}}
\date{}
\maketitle 
\begin{abstract}
We study the algebraic aspects of equivariant quantum cohomology algebra 
of the flag manifold. We introduce and study the quantum double
Schubert polynomials $\wt{\s}_w(x,y)$, which are the 
Lascoux--Sch\"utzenberger type representatives of the equivariant
quantum cohomology classes. Our approach is based on the quantum Cauchy 
identity. We define also quantum Schubert polynomials $\wt{\s}_w(x)$ as
the Gram--Schmidt orthogonalization of some set of monomials with respect
to the scalar product, defined by the Grothendieck residue. Using quantum
Cauchy identity, we prove that $\wt{\s}_w(x)=\wt{\s}_w(x,y)|_{y=0}$ and
as corollary obtain a simple formula for the quantum Schubert polynomials
$\wt{\s}_w(x)=\partial_{ww_0}^{(y)}\wt{\s}_{w_0}(x,y)|_{y=0}$\,.
We also prove the higher genus analog of Vafa--Intriligator's formula for 
the flag manifolds and study the quantum residues generating function. 
We introduce the extended Ehresman--Bruhat order on the symmetric 
group and formulate the equivariant quantum Pieri rule.
\end{abstract}
\newpage

\section{Introduction.}

The structure constants of the quantum cohomology ring 
are given by the third derivatives of the Gromov-Witten 
potential $F.$ The Gromov-Witten potential $F$ is a 
generating function of the Gromov-Witten invariants. 
The axioms of the tree level Gromov-Witten invariants 
$$ \langle I_{0, m, \beta}^{V} \rangle : 
H^{\ast}(V,{\bf Q})^{\otimes m} \longrightarrow {\bf Q}, 
$$
$\beta \in H_2 (V,{\bf Z}),$ 
for a target space $V$ are given by Kontsevich and Manin [KM]. 
Let $X_1,\ldots ,X_m$ be cycles on $V$ and $X_1^{\ast},\ldots ,
X_m^{\ast}$ their dual classes. Then the invariant 
$\langle I_{0,m,\beta}^V \rangle (X_1^{\ast}\otimes \cdots 
\otimes X_m^{\ast} )$ can be considered as the virtual 
number of the stable maps $f$ from $m$-pointed rational curve 
$({\bf P}^1 ; p_1,\ldots ,p_m)$ to $V,$ such that 
the image of $f$ represents the homology class $\beta$ and 
$f(p_i)\in X_i.$ 

In case of flag variety $Fl_n:=SL_n/B$ of type $A_{n-1}$ the 
potential $F$ is given as follows. 
Let $\Omega _v$ be the dual class of the Schubert cycle $X_v$ 
corresponding to a permutation $v\in S_n.$ Then the 
potential $F_{\omega}((t_v)_{v\in S_n})$ is defined by 
$$ F_{\omega}((t_v)_{v\in S_n}) = \sum_{\beta} 
\sum_{m=\Sigma m_v \geq 3} \exp(-\int_{\beta}\omega)
\frac{ \langle I_{0,m,\beta}^{V} \rangle 
(\ds\bigotimes _{v\in S_n}\Omega_v^{\otimes m_v})}
{\ds\prod_{v\in S_n} m_v !}\prod_{v\in S_n}t^{m_v}, 
$$
where $\omega$ is a K\"ahler form. For each point 
$t\in H^{\ast} (Fl_n),$ the quantum multiplication law is 
given by 
$$ \Omega_u \ast \Omega_v = \sum_{w\in S_n}
\frac{\partial ^3 F_{\omega}}{\partial t_u \partial t_v 
\partial t_w}(t) \Omega_{ww_0}, 
$$ 
where $w_0$ is the permutation of maximal length. 
The algebra with this multiplication law is called a quantum cohomology 
ring, 
which is denoted by $QH^{\ast}_t (Fl_n).$ 
The associativity  of the quantum multiplication is 
equivalent to the Witten--Dijkgraaf--Verlinde--Verlinde 
(WDVV) equation 
\[ \sum_{v\in S_n}
\frac{\partial ^3 F_{\omega}}
{\partial t_{u_1} \partial t_{u_2} \partial t_{v}}
\frac{\partial ^3 F_{\omega}}
{\partial t_{vw_0} \partial t_{u_3} \partial t_{u_4}}
= \sum_{v\in S_n}
\frac{\partial ^3 F_{\omega}}
{\partial t_{u_2} \partial t_{u_3} \partial t_{v}}
\frac{\partial ^3 F_{\omega}}
{\partial t_{vw_0} \partial t_{u_1} \partial t_{u_4}}, \]
for any $u_1,u_2,u_3,u_4 \in S_n.$ 
From [KM, Proposition 4.4], the potential $F_{\omega}$ 
satisfies 
\[ F_{\omega} (t) = F_{\omega - t^{(2)}}(t-t^{(2)}), \] 
where $t^{(2)} = \ds\sum_{l(v) = 1} t_v \Omega_v .$ 
Hence, we may assume $\omega = 0.$ 
The potential $F$ is decomposed as a sum of the classical 
part $f_{\rm cl}$ and the quantum correction $f,$ where 
\[ f_{\rm cl} = \frac{1}{6} t_{id} 
\left( \sum_{u,v}\delta_{u,vw_0}t_u t_v \right) , \] 
and 
\[f = \sum_{\beta}\sum_{m = \Sigma m_v \geq 3}
\frac{\langle I_{0,m,\beta}^{Fl} \rangle 
(\ds\bigotimes _{l(v) \geq 1}\Omega _v^{\otimes m_v})}
{\ds\prod_{l(v)\geq 1} m_v !}\prod_{l(v)\geq 1}t^{m_v}. \]
From the axioms of the Gromov-Witten invariants
([KM,(2.2.4)]), we have 
\[ \langle I_{0,m,\beta}^{Fl} \rangle 
(\bigotimes _{l(v)\geq 1} 
\Omega _{v}^{\otimes m_v}) = \langle 
I_{0,m-\Sigma_{l(u)=1}m_u,\beta}^{Fl} 
\rangle (\bigotimes _{l(v) > 1} \Omega _{v}^{\otimes m_v} )
\prod_{l(u)=1} \left( \int_{\beta}\Omega _u \right). \]
Hence, the quantum correction $f$ is expressed as 
\[ f = \sum_{m_v, b_u} N((m_v)_{l(v)>1} \mid (b_u)_{l(u)= 1})
\exp(\sum_{l(u)=1}(b_u t_u))\prod_{l(v)>1}\frac{t_v^{m_v}}
{m_v !}, \]
where \[ N((m_v) \mid (b_u)) = \langle 
I_{0, \Sigma m_v, \Sigma b_u X_u}^{Fl} \rangle 
\left( \bigotimes _{l(v)>1} \Omega _v^{\otimes m_v} \right). \]
If $n=3,$ the WDVV equations with the initial condition 
\[ N(0,0,1 \mid 1,0) = N(0,0,1 \mid 0,1) = 1 \] 
or 
\[ N(2,0,0 \mid 1,0) = N(0,2,0 \mid 0,1) = 1 \] 
determine all the coefficients $N(\lambda, \mu, \nu \mid a,b)$ 
uniquely. In fact, Di Francesco and Itzykson [FI] gave the 
coefficients $N(\lambda, \mu, \nu \mid a,b)$ for $a+b \leq 10.$ 

Let $x_1 = \Omega _{(1,2)}$ and $x_i = \Omega _{(i,i+1)} - 
\Omega _{(i-1,i)}$ for $2\leq i\leq n-1,$ where $(i,j)$ is 
the transposition that interchanges $i$ and $j.$ 
Since the classical cohomology ring is generated by 
$x_1,\ldots ,x_{n-1},$ 
the quantum cohomology ring $QH^{\ast}_t(Fl_n)$ 
is also generated by $x_1,\ldots ,x_{n-1}$ in 
the neighborhood of the origin $t=0.$ Moreover, 
we can choose a neighborhood of the origin on 
which $QH^{\ast}_t(Fl_n)$ is a complete 
intersection ring. 
Then, let ${\cal I}\subset {\bf C}[x_1,\ldots 
,x_{n-1}]$ be the defining ideal of the quantum cohomology 
ring $QH^{\ast}_t(Fl_n).$ The Schubert class $\Omega _v$ is 
expressed by the Schubert polynomial ${\s}_v(x_1,\ldots ,
x_{n-1})$ in the classical cohomology ring. However, in the 
quantum cohomology ring, the class corresponding to the 
Schubert polynomial ${\s}_v$ is no longer the Schubert 
class $\Omega _v.$ Hence, the polynomial $\wt{\s}_{v}
^{t}(x_1,\ldots ,x_{n-1})$ expressing $\Omega _v$ gives a 
deformation of the Schubert polynomial. We call it a big
quantum Schubert polynomial. We identify the residue 
pairing defined by ${\cal I}$ with the intersection form on the 
cohomology ring, so the big quantum Schubert polynomials are 
obtained by orthogonalization the basis consisting of the 
classical Schubert polynomials. 

It is difficult to describe the defining ideal $\cal I$ of the big
quantum cohomology ring for generic $t,$ so the big quantum 
Schubert polynomials are complicated in general. However, 
in the case where parameters $t_v = 0$ for all permutations 
$v\in S_n$ such that 
$l(v) > 1,$ the defining relations of the quantum cohomology 
ring are known by the results of A. Givental and B. Kim [GK] and 
I. Ciocan-Fontanine [C]. We call it the small quantum cohomology 
ring. The structure constants of the small quantum cohomology 
ring are given by 
\[ 
\frac{\partial ^3 F}{\partial t_{v_1} \partial t_{v_2} 
\partial t_{v_3}} (t^{(2)})= \]  
\begin{eqnarray*}
& = &
\sum_{\beta}\sum_{m_u \geq 0}
\frac{\langle I_{0,\Sigma m_u +3,\beta}^{Fl}
\rangle (\Omega _{v_1}
\otimes \Omega _{v_2} \otimes \Omega _{v_3} \otimes \ds(\bigotimes 
_{l(u) =1}\Omega_{u}^{\otimes m_u}))}{\ds\prod_{l(u)=1} m_u !}
\prod_{l(u)=1}t_u^{m_u} \\ 
& = & \sum_{\beta} \langle I_{0,3,\beta}^{Fl} \rangle 
(\Omega_{v_1}\otimes \Omega_{v_2}\otimes \Omega_{v_2}) 
e^{\Sigma b_i t_{(i,i+1)}}, 
\end{eqnarray*} 
where the sum runs over 
\[ \beta = b_1 X_{(1,2)}+\cdots + b_{n-1}X_{(n-1,n)} \] 
with $b_i\in {\bf Z}_{\geq 0}.$ Hence the small quantum cohomology ring is 
determined by the invariants 
$\langle I_{0,3,\beta}^{Fl} \rangle .$ 
For a monomial $x_{i_1}\ast \cdots \ast x_{i_m},$ the 
$m$-point correlation function determined by 
the small quantum cohomology ring (the so-called small quantum cohomology 
ring correlation function) is defined to be
\[ \langle x_{i_1} \cdots x_{i_m} \rangle 
= \int_{Fl_n}x_{i_1}\ast \cdots \ast x_{i_m}. \] 
 If $m\geq 4,$ the small quantum cohomology ring correlation function can
be expressed as follows: 
\[ \langle x_{i_1} \cdots x_{i_m} \rangle = \] 
\[ \sum_{\beta = \beta_1 + \cdots + \beta_{m-3}} 
\sum_{v_1,\ldots ,v_{m-3}} 
e^{ -\int_{\beta}\omega } 
\langle I_{0,3,\beta_1} \rangle (x_{i_1} \otimes x_{i_2} 
\otimes \Omega_{v_1} ) 
\langle I_{0,3,\beta_2} \rangle (\Omega_{v_1 w_0} 
\otimes x_{i_3} \otimes \Omega_{v_2} ) \] 
\[ \cdots 
\langle I_{0,3,\beta_{m-2}} \rangle (\Omega_{v_{m-3} w_0} 
\otimes x_{i_{m-1}} \otimes x_{i_m} ). \] 
Hence, if $m\ge 4$, the Gromov-Witten invariants 
$\langle I_{0,m,\beta} \rangle 
(x_{i_1} \otimes \cdots \otimes x_{i_m})$
do not appear as coefficients of small quantum cohomology ring 
correlation functions. 

Let us explain briefly the main results obtained in our paper. 
Follow to A.~Givental and B.~Kim [GK], and 
I.~Ciocan--Fontanine [C], we define the 
quantum elementary symmetric polynomials 
$\wt e_1,\ldots ,\wt e_n$ by the formula 
$$ \det\left(\begin{array}{ccccccc} x_1+t & q_1&0 &\ldots &\ldots 
&\ldots &0\\
-1 & x_2+t & q_2 &0 &\ldots &\ldots & 0\\
0 & -1 &x_3+t & q_3 & 0 &\ldots & 0 \\
\vdots &\ddots &\ddots & \ddots &\ddots &\ddots &\vdots \\
0&\ldots & 0 &-1&x_{n-2}+t &q_{n-2} & 0 \\
0 &\ldots &\ldots & 0 &-1 & x_{n-1}+t & q_{n-1}\\
0 &  \ldots &\ldots &\ldots & 0 & -1 & x_n+t
\end{array}\right) 
$$ 
$$= t^n + \wt e_1t^{n-1} + \wt e_2t^{n-2}+\cdots +\wt e_{n},
$$  
where $q_i = e^{t_{(i,i+1)}}.$ The defining ideal $\wt I$ of the 
small quantum cohomology ring is generated by the 
quantum elementary symmetric polynomials, namely 
{\begin{equation}QH^*(Fl_n,{\bf Z}):=QH^{\ast}(Fl_n) = {\bf Z}[x_1,\ldots ,x_n;q_1,
\ldots ,q_{n-1}]/(\wt e_1,\ldots ,\wt e_n). \label{211}
\end{equation}

In the classical case $q_1=\cdots =q_{n-1}=0$, on the quotient ring
$$\a :={\bf Z}[ x_1,\ldots ,x_n]/(e_1(x),\ldots ,e_n(x))\simeq 
H^*(Fl_n,{\bf Z})
$$
there exists a natural pairing $\langle f,g\rangle =\eta 
(\partial_{w_0}(fg))$ which comes from the intersection pairing in the
homology group $H_*(Fl_n,{\bf Z})$ of the flag variety. We can interpret
the pairing $\langle ,\rangle$ as the Grothendieck residue pairing with 
respect to the ideal $I$ (see Subsection 2.5):
$$\langle f,g\rangle =\Res_I(fg),
$$
where $I=\wt I|_{q=0}$.

Our first observation is that a natural residue pairing (we call it the
quantum residue pairing)
$$\langle f.g\rangle_Q=\Res_{\wt I}(fg)
$$
on the quotient ring $\overline{\a} :={\bf Z}[x_1,\ldots ,x_n]/\wt I$ 
corresponds to the intersection pairing in quantum cohomology 
$QH^*(Fl_n,{\bf Z})$ under a natural isomorphism (\ref{211}).

It is well-known (e.g. [LS2], [M]) that the classical Schubert polynomials
form an orthonormal basis (with respect to the pairing
$\langle ,\rangle$) in the cohomology ring of flag manifold and also
give a linear basis in the quantum cohomology ring $QH^*(Fl_n,{\bf Z})$,
[GK], [MS].
However, the classical Schubert polynomials do not orthogonal with respect
to the quantum pairing any more. Thus, it is natural to ask: what kind of
polynomials one can obtain applying the Gram--Schmidt orthogonalization 
to the classical Schubert polynomials with respect to the quantum pairing
$\langle ,\rangle_Q$? Omiting some details with ordering (see Definition 5),
the answer is: quantum Schubert polynomials.

Our second observation is: to work with the equivariant quantum
cohomology algebra ([GK], [K2]) is more convenient then with quantum 
cohomology ring itself.
The main reason is that one can 
find Lascoux--Sch\"utzenberger's type representative for any equivariant 
quantum cohomology class. In other words, each quantum double Schubert 
polynomial $\wt{\s}_w(x,y)$ can be obtained from the top one by using the 
divided difference operators acting on the $y$ variables. \medskip \\
{\bf Theorem-Definition A} \ {\it Let $x=(x_1,\ldots ,x_n)$, 
$y=(y_1,\ldots ,y_n)$ 
be two sets of variables, and
$$\wt{\s}_{w_0}^{(q)}(x,y):=\prod_{i=1}^{n-1}\Delta_i(y_{n-i}~|~x_1,\ldots 
,x_i),
$$
where $\Delta_k(t~|~x_1,\ldots 
,x_k):=\ds\sum_{j=0}^kt^{k-j}e_j(x_1,\ldots ,x_k~|~q_1,\ldots ,q_{k-1})$
is the generating function for the quantum elementary symmetric functions 
in $x_1,\ldots ,x_k$.

Then \ $\wt{\s}_w^{(q)}(x,y)=\partial_{ww_0}^{(y)}\wt{\s}_{w_0}^{(q)}(x,y)$.}

We define the quantum Schubert polynomials $\wt{\s}_w^{(q)}(x)$ as 
Gram--Schmidt's orthogonalization of the set of lexicographically ordered 
monomials\break
$\{ x^I \mid I\subset (n-1,n-2,\ldots ,1,0) \}$
with respect to the quantum pairing $\langle ,\rangle_Q$, see 
Definition 5. One of our main results is the quantum analog of Cauchy's 
identity for (classical) Schubert polynomials, [M], (5.10). \medskip \\
{\bf Theorem B} \ {\it (Quantum Cauchy's identity)
\begin{equation}\sum_{w\in S_n}\wt{\s}_w^{(q)}(x)\s_{ww_0}(y)=
\wt{\s}_{w_0}^{(q)}(x,y).
\end{equation}}

We give a geometric proof of Theorem B in Section 7 using the arguments 
due to I.~Ciocan-Fontanine [C]; more particularly, we reduce directly a 
proof of Theorem B to that of the following geometric statement: \medskip \\
{\bf Lemma} \ {\it Let $I\subset\delta =(n-1,n-2,\ldots ,1,0)$ and $w\in 
S_n$ be a permutation, then
\begin{equation} \langle\wt e_I(x),\wt{\s}_w(x)\rangle_Q=\langle 
e_I(x),\s_w(x)\rangle ,
\end{equation}
where $e_I(x):=\ds\prod_{k=1}^{n-1}e_{i_k}(x_1,\ldots ,x_{n-k})$\\
(resp. 
$\wt e_I(x):=\ds\prod_{k=1}^{n-1}\wt e_{i_k}(x_1,\ldots ,x_{n-k}~|~q_1,\ldots 
,q_{n-k-1})$)\\ 
is the elementary polynomial (resp. quantum elementary 
polynomial), see Section 5.2.}

It is the formula (3) that we prove in Section 7 using the geometrical 
arguments from [C] and [K2]. By product, it follows from our proof that 
quantum Schubert polynomials $\hat{\s}_w(x)$ defined geometrically (see 
Section~6) coincide with those defined algebraically (see Definition 5):
$$\hat{\s}_w(x) \equiv \wt{\s}_{w^{-1}}(x) \ (\mod \wt I).
$$
It is interesting to note, that the intersection numbers $\langle 
e_I(x),\s_w(x)\rangle$ (which are nonnegative!) are precisely the 
coefficients of corresponding Schubert polynomial:
$$\s_w(x)=\sum_{I\subset\delta}\langle e_I(x),\s_w(x)\rangle x^{\delta 
-I}.
$$
The quantum Cauchy formula (2) plays the important role in our approach 
to the quantum Schubert polynomials. As a direct consequence of (2), we 
obtain the Lascoux--Sch\"utzenberger type formula for quantum Schubert 
polynomials (cf. Theorem-Definition A). \medskip \\
{\bf Theorem C} \ {\it Let $\wt{\s}_{w_0}(x,y)$ be as in 
Theorem-Definition A, then}
$$\wt{\s}_{w}(x)=\partial_{ww_0}^{(y)}\wt{\s}_{w_0}(x,y)|_{y=0}.
$$

In Section 5 we introduce a quantization map 
\[ P_n \rightarrow \overline{P}_n, \; \; f \mapsto \wt{f}. \] 
The quantization is a linear map which preserves the pairings, i.e., 
\[ \langle \wt{f},\wt{g} \rangle _{Q} = \langle f,g \rangle , \; \; 
f,g \in P_n. \] 
Using the quantum Cauchy formula (2), we prove that quantum double 
Schubert polynomials are the quantization of classical ones. 
Another class of polynomials having a nice quantization is the set of 
elementary polynomials 
\[ e_I(x) := \prod_{k=1}^{n-1}e_{i_k}(x_1,\ldots x_{n-k}), \; \; 
I=(i_1,\ldots ,i_{n-1})\subset \delta . \] 
It follows from Theorem B that quantization $\wt{e}_I(x)$ of elementary 
polynomial $e_I(x)$ is given by 
\[ \wt{e}_I(x) = \prod_{k=1}^{n-1} e_{i_k}(x_1,\ldots ,x_{n-k}
\mid q_1,\ldots ,q_{n-k-1} ). \] 
More generally, we make a conjecture (''quantum Schur functions '') that 
quantization of the flagged Schur function (see [M], (3.1), (4.9) and (6.16)) 
\[ s_{\lambda / \mu }(X_1,\ldots ,X_n) = 
{\rm det}\left( h_{\lambda_i-\mu_j-i+j}(X_i)\right) _{1\leq i,j \leq n} \] 
is given by 
\[ \wt{s}_{\lambda / \mu }(X_1,\ldots ,X_n) = 
{\rm det}\left( \wt{h}_{\lambda_i-\mu_j-i+j}(X_i)\right) 
_{1\leq i,j \leq n}, \] 
where $\wt{h}_k(X)$ is the quantum complete homogeneous symmetric function 
of degree $k,$ and $X_1\subset \cdots \subset X_n$ are the flagged sets of 
variables (see Section 5). 

In Section 5.2 we consider a problem how to quantize monomials. It seems 
to be difficult to find an explicit determinantal formula for a quantum 
monomial $\wt{x}^I,$ i.e., to find a quantum analog of the 
Billey-Jockusch-Stanley formula for Schubert polynomials in terms of 
compatible sequences [BJS]. We prove the following formulae for 
quantum monomials 
\[ \wt{x}^I = \sum_{w\in S_n} \eta (\partial _{w}x^I ) \wt{\s}_w(x), 
\; \; I\subset \delta, \] 
\[ \wt{\s}_{w_0}(x,y) = \sum_{I\subset \delta}\wt{x}^I 
e_{\delta - I}(y). \] 

In section 8.1 we give a proof of the higher genus analog of the Vafa--\break
Intriligator type formula for the flag manifold. 

In Section 8.3 we study a problem how to compute the quantum residues. 
This is important for computation of small quantum cohomology ring 
correlation functions (or correlation functions, for short) and the
Gromov--Witten invariants, see Introduction and Theorem 11. We introduce 
the generating function 
\[ \Psi(t)= \langle \prod_{i=1}^{n-1}\frac{t_i}{t_i - x_i} \rangle \]
for quantum residues and give a characterization of this function as the 
unique solution to some system of differential equations, see Proposition 14. 
In Appendix~B we calculate the generating function $\Psi(t)$ for 
the case $n=3$ explicitly. 

In Section 9 we introduce the extended Ehresman--Bruhat order and 
give a sketch of a proof of equivariant quantum Pieri rule. Details 
will appear elsewhere. 

In Appendix A one can find a list of explicit expressions for the quantum
double Schubert polynomials for the symmetric group $S_4$.

We would like to mention, that in the recent preprint 
``Quantum Schubert polynomials'' by
S.~Fomin, S.~Gelfand and A.~Postnikov, [FGP], developed a different
approach to the theory of quantum Schubert polynomials, based on 
the remarkable family of commuting operators $X_i$ ([FGP], (3.2)). 
Among main results, obtained by S.~Fomin, S.~Gelfand and A.~Postnikov, are
definitions, orthogonality, quantum Monk's formula and other 
properties of quantum Schubert polynomials; definition of quantization map 
and quantum multiplication. 

Besides some overlap with the preprint of S.~Fomin, S.~Gelfand
and A.~Postnikov, our works were done independently and based on 
the different approaches, which allow to obtain the mutually complementary
results. 

\section{Classical Schubert polynomials.}

In this section we give a brief review of the theory of Schubert 
polynomials created by A.~Lascoux and M.-P.~Sch\"utzenberger. In 	
exposition we follow to the I.~Macdonald book [M1] where proofs and more 
details can be found.

\subsection{Divided differences.}

Let $x_1,\ldots ,x_n,\ldots$ be independent variables, and let
$$P_n:={\bf Z}[x_1,\ldots ,x_n]
$$
for each $n\ge 1$, and 
\begin{equation}
P_{\infty}:={\bf Z}[x_1,x_2,\ldots ]={\ds\bigcup_{n=1}^{\infty}P_n}.
\end{equation}
Let us denote by $\Lambda_n:={\bf Z}[x_1,\ldots ,x_n]^{S_n}\subset P_n$
the ring of symmetric polynomials in $x_1,\ldots ,x_n$, and by 
$H_n:=\{\ds\sum_{I=(i_1,\ldots ,i_n)}a_Ix^I~|~a_I\in{\bf Z}, \ 0\le 
i_k\le n-k, \forall k\}$ the additive subgroup of $P_n$ spanned by all 
monomials $x^I:=x_1^{i_1}x_2^{i_2}\ldots x_n^{i_n}$ with $I\subset\delta 
:=\delta_n=(n-1,n-2,\ldots ,1,0)$. For 
$1\le i\le n-1$ let us define a linear operator $\partial_i$ acting on 
$P_n$
\begin{equation}
(\partial_i f)(x)= \frac{f(x_1,\ldots,x_i,x_{i+1},\ldots,x_n)
-f(x_1,\ldots,x_{i+1},x_i,\ldots,x_n)}
{x_i-x_{i+1}}.
\end{equation} 
Divided difference operators $\partial_i$ satisfy the following 
relations 
\begin{eqnarray}
\partial_i^2 & = & 0,\nonumber \\ 
\partial_i \partial_j & = & \partial_j \partial_i,
\ \ {\rm if} \mid i-j \mid >1, \\ 
\partial_i \partial_{i+1} \partial_i & = & 
\partial_{i+1} \partial_i \partial_{i+1}, \nonumber
\end{eqnarray}
and the Leibnitz rule 
\begin{equation} 
\partial_i (fg) = \partial_i(f)g+s_i(f) \partial_i(g).\label{5}
\end{equation}
It follows from (\ref{5})  that $\partial_i$ is a $\Lambda_n$-linear operator. 

For any permutation $w\in S_n,$ let us denote by $R(w)$ the set of 
reduced words for $w,$ i.e. sequences $(a_1,\ldots,a_p)$ such that 
$w= s_{a_1}\cdots s_{a_p},$ where $p=l(w)$ is the length of permutation 
$w\in S_n$, and $s_i=(i,i+1)$ is the simple 
transposition that interchanges $i$ and $i+1.$ 

For any sequence ${\bf a}=(a_1,\ldots,a_p)$ of positive integers, 
we define 
\[ \partial_{\bf a}= \partial_{a_1}\cdots\partial_{a_p}. \] 
\begin{pr}\label{2} ([M1], (2.5),(2.6))

$\bullet$ If ${\bf a,b}\in R(w),$ then $\partial_{\bf a}=\partial_{\bf b}.$ 

\vskip 0.1cm
$\bullet$ If ${\bf a}$ is not reduced, then $\partial_{\bf a}=0.$ 
\end{pr}
From Proposition \ref{2} it follows that an operator \[ \partial_w = \partial_{\bf a} \] 
is well-defined, where ${\bf a}$ is any reduced word for $w.$ By (\ref{5}), the 
operators $\partial_w$, $w\in S_n$, are $\Lambda_n$ linear, i.e. if
$f\in\Lambda_n$, then
$$\partial_w(fg)=f\partial_w(g).
$$                                     

\subsection{Schubert polynomials.}

Let $\delta =\delta_n=(n-1,n-2,\ldots ,1,0)$, so that 
$x^{\delta}=x_1^{n-1}x_2^{n-2}\ldots x_{n-1}$.

\begin{de} (Lascoux--Sch\"utzenberger [LS1]).
 For each permutation \hbox{$w\in S_n$} the Schubert polynomial 
 ${\s}_w$ is defined to be 
 $${\s}_w(x)=\partial_{w^{-1}w_0}(x^{\delta}),
 $$
where $w_0$ is the longest element of $S_n$.
 \end{de}

\begin{pr} ([M1], (4.2),(4.5),(4.11),(4.15)).\\
$\bullet$ Let $v,w\in S_n$. Then 
$$\partial_v\s_w=\left\{\begin{array}{ll}
                        \s_{wv^{-1}}, & \mbox{\rm if \ $l(wv^{-1})=l(w)-l(v)$,}\\
                        0, & \mbox{\rm otherwise}.
\end{array}\right.
$$
$\bullet$ (Stability). Let $m>n$ and let $i:S_n\hookrightarrow S_m$ to be the 
natural embedding. Then
$$\s_w=\s_{i(w)}.
$$
$\bullet$ The Schubert polynomials $\s_w, w\in S_n$ form a ${\bf 
Z}$--basis of $H_n$.\\
$\bullet$ (Monk's formula). Let $f=\ds\sum_{i=1}^n\alpha_ix_i$, 
$w\in S_n$. Then
\begin{eqnarray*}f\s_w &=&\sum (\alpha_i-\alpha_j)\s_{wt_{ij}},\\
\partial_w(fg) &=& w(f)\partial_wg+\sum (\al_i-\al_j)\partial_{wt_{ij}}g,
\end{eqnarray*}
where $t_{ij}$ is the transposition that interchanges $i$ and $j$, and both 
sums are over all pairs $i<j$ such that $l(wt_{ij})=l(w)+1$.
\end{pr}

\subsection{Scalar product.}

Let us define a scalar product on $P_n$ with values in $\Lambda_n$, by 
the rule
\begin{equation}
\langle f,g\rangle =\partial_{w_0}(fg),~f,g\in P_n,\label{100}
\end{equation}
where $w_0$ is the longest element of $S_n$.

The scalar product $\langle ,\rangle$ defines a non-degenerate pairing
$\langle ,\rangle_0$ on the quotient ring $P_n/I_n\cong H^*(Fl_n,{\bf Z})$,
where $I_n$ is the ideal in $P_n$ generated by the elementary symmetric
polynomials $e_1(x),\ldots ,e_n(x)$.

\begin{pr} ([M1], (5.3),(5.4),(5.6),(4.13),(5.10)).\\ ~\\
$\bullet$ If $f\in\Lambda_n$, then $\langle fh,g\rangle =f\langle 
h,g\rangle $;\\ 
$\bullet$ If $f,g\in P_n$, $w\in S_n$, then 
$\langle\partial_wf,g\rangle =\langle f,\partial_{w^{-1}}g\rangle$;\\ ~\\
$\bullet$ (Orthogonality) If $l(u)+l(v)=l(w_0)$, then 
$\langle\s_u,\s_v\rangle =\left\{\begin{array}{ll}
                                        1, & \mbox{\rm if \ $u=w_0v$,}\\
                                        0, & \mbox{\rm otherwise}.
\end{array}\right.$\\
$\bullet$ The Schubert polynomials $\s_w$, $w\in S_n$, form a 
$\Lambda_n$--basis of $P_n$;\\ ~\\
$\bullet$ The Schubert polynomials $\s_w$, $w\in S^{(n)}$, form a ${\bf 
Z}$--basis of $P_n$, where for each $n\ge 1$, $S^{(n)}$ is the set of all 
permutations $w$ such that the code $w$ has length $\le n$;\\
$\bullet$ (Cauchy's formula)
$$\sum_{w\in S_n}\s_w(x)\s_{ww_0}(y)=\prod_{i+j\le n}(x_i+y_j).
$$
\end{pr}
\begin{pr} Schubert polynomials are uniquely characterized by the 
following properties\\
1. (Orthogonality) $\langle\s_u,\s_v\rangle_0 =\left\{\begin{array}{ll}
                                        1, & \mbox{\rm if \ $u=w_0v$,}\\
                                        0, & \mbox{\rm otherwise}.
\end{array}\right.$\\
2. Let $w$ be a permutation in $S_n$ and $c(w)=(c_1,c_2,\ldots ,c_n)$ its 
code, then
$$\s_w(x)=x^{c(w)}+\sum\al_Ix^I,
$$
where $I\subset\delta$, $\al_I>0$ and $I$ lexicographically smaller then 
$c(w)$.
\end{pr} 
{\bf Remark 1} \ 1) (Definition of the code, [M1], p.9). \\ 
For a permutation $w\in S_n,$ we define 
\[ c_i = \sharp \{ j \mid i<j, w(i)> w(j) \} . \] 
The sequence $c(w)=(c_1,c_2,\ldots ,c_n)$ is called the 
code of $w.$ \smallskip \\ 
2) Schubert polynomials are obtained as 
Gram-Schmidt's orthogonalization of the set of monomials 
~$\{ x^I\}_{I\subset \delta}$~ ordered lexicographically. 

\subsection{Double Schubert polynomials.}

Let $x=(x_1,\ldots ,x_n)$, $y=(y_1,\ldots ,y_n)$ be two sets of 
independent variables, and
$$\s_{w_0}(x,y):=\prod_{i+j\le n}(x_i+y_j).
$$

\begin{de} (Lascoux--Sch\"utzenberger [LS2]).
For each permutation \hbox{$w\in S_n$,} the double Schubert polynomial 
$\s_w(x,y)$ is defined to be
$$\s_w(x,y)=\partial_{w^{-1}w_0}^{(x)}\s_{w_0}(x,y),
$$
where divided difference operator $\partial_{w^{-1}w_0}^{(x)}$ 
acts on the $x$ variables.
\end{de}

\begin{pr} ([M1], (6.3),(6.8)).\\
$\bullet$ $\s_w(x,y)=\ds\sum_u\s_u(x)\s_{uw^{-1}}(y)$, summed over all $u\in 
S_n$, such that\break $l(u)+l(uw^{-1})=l(w)$;\\
$\bullet$ (Interpolation formula). For all $f\in{\bf Z}[x_1,\ldots ,x_n]$ 
we have
$$f(x)=\sum_w\s_w(x,-y)\partial_w^{(y)}f(y)
$$
summed over all permutations $w\in S^{(n)}$.
\end{pr}

Double Schubert polynomials appear in algebra and geometry as cohomology 
classes related to degeneracy loci of flagged vector bundles. If $h:E\to 
F$ is a map of rank $n$ vector bundles on a smooth variety $X$,
$$E_1\subset E_2\subset\cdots\subset E_n=E, \ F:=F_n\to 
F_{n-1}\to \cdots \to F_1
$$
are flags of subbundles and quotient bundles, then there is a degeneracy 
locus $\Omega_w(h)$ for each permutation $w$ in the symmetric group 
$S_n$, described by the conditions
$$\Omega_w(h)=\{ x\in X~|~{\rm rank}(E_p(x)\to F_q(x))\le\#\{ i\le q, 
w_i\le p\} , \forall p,q\}.
$$
For generic $h$,\  $\Omega_w(h)$ is irreducible, codim $\Omega_w(h)=l(w)$, 
and the class $[\Omega_w(h)]$ of this locus in the Chow ring of $X$ is 
equal to the double Schubert polynomial $\s_{w_0w}(x,-y)$, where
\begin{eqnarray}  x_i &=& c_1(\ker (F_i\to F_{i-1})),\nonumber\\
y_i &=& c_1(E_i/E_{i-1}), \ 1\le i\le n.\nonumber
\end{eqnarray}
It is well-known [F] that the Chow ring of flag variety $Fl_n$ admits the 
following description
$$CH^*(Fl_n)\cong{\bf Z}[x_1,\ldots ,x_n,y_1,\ldots ,y_n]/J,
$$
where $J$ is the ideal generated by
$$e_i(x_1,\ldots ,x_n)-e_i(y_1,\ldots ,y_n), \ 1\le i\le n,
$$
and $e_i(x)$ is the $i$-th elementary symmetric function in the variables 
$x_1,\ldots ,x_n$.

$\bullet$ ([LS2], [KV]) The ring ${\bf Z}[x_1,\ldots ,x_n,y_1,\ldots ,y_n]/J$ is 
a free module of dimension $n!$ over the ring $R$, with basis either 
$\s_w(x)$, or $\s_w(x,y),\ w\in S_n$, where
$$R:=\frac{{\bf Z}[x_1,\ldots ,x_n]\otimes{\bf Sym}[y_1,\ldots ,y_n]}
{J}.
$$

\subsection{Residue pairing.}

Let $I$ be an ideal in $\bar P_n =R[x_1,\ldots ,x_n]$, $R\subset {\bf C}$, 
generated by a regular system of parameters $\varphi_1,\ldots ,\varphi_n$, 
and $\a :=\bar P_n /I$.

\begin{pr} ([GH], [EL]).

$\bullet$ $\dim_R\a <\infty$.

$\bullet$ ${\cal H}:=\det \left(\ds\frac{\partial\varphi_i}
{\partial x_j}\right)\not\in I$.
\end{pr}

Let $d_0:=\deg\cal H$, where we assume that $\deg x_i=1$ for all 
$1\le i\le n$.

\begin{pr} ([EL])

$\bullet$ If $f\in\bar P_n$ and $\deg f=d_0$, then there exists a non 
zero $\al\in R$ such that
$$f\equiv\frac{\al}{n!}{\cal H} \ (\mod I).
$$

$\bullet$ If $f\in\bar P_n$, $f\ne 0$, and $\deg f>d_0$, then there exists 
$g\in\bar P_n$ such that $\deg g\le d_0$ and $g\equiv f \ (\mod I)$.
\end{pr}

\begin{de} (Grothendieck residue with respect to the ideal $I$).

Let $f\in\bar P_n$ and $\deg f<d_0$, then we define
$$\Res_I(f)=0.
$$

If $\deg f=d_0$, then $f\equiv \ds\frac{\al}{n!}{\cal H} \ (\mod I)$ and we 
define $\Res_I(f):=\al$.

Finally, if $\deg f>d_0$, then choose $g\in\bar P_n$ such that $g\equiv 
f \ (\mod I)$ and $\deg g\le d_0$, and define
$$\Res_I(f):=\Res_I(g).
$$
\end{de}

We will use also notation $\langle f\rangle_I$ instead of $\Res_I(f)$.

Finally, let us define a residue pairing $\langle , \rangle_I$ on 
$\overline P_n$ using the Grothendieck residue
$$\langle f,g\rangle_I=\Res_I(f,g), \ \ f,g\in\overline P_n.
$$
\begin{pr} ([GH]).\\
$\bullet$ If $f\in I$, then $\Res_I(f)=0$.\\
$\bullet$ The residue pairing $\langle , \rangle_I$ induces a 
non-degenerate pairing on $\a =\overline P /I$.
\end{pr}

We will use this general construction of residue pairing in the following 
two cases:

i) $R={\bf Z}$, $I_n\subset P_n$ is an ideal generated by elementary 
symmetric polynomials $e_1(x), \ldots ,e_n(x)$. It is well-known that if 
$\Fl :=SL(n)/B$ is the flag variety of type $A_{n-1}$, then
$$H^*(\Fl ,{\bf Z})\simeq P_n/I_n,
$$
and residue pairing $\langle ,\rangle$ on $P_n/I_n$ coincides with the 
scalar product on $P_n/I_n$ induced by (\ref{100}).

ii) $R={\bf Z}[q_1,\ldots ,q_{n-1}]$, $\wt I_n\subset\overline P_n$ 
is an ideal 
generated by the quantum elementary symmetric functions $\wt 
e_1(x),\ldots ,\wt e_n(x)$. It is a result of A.~Grivental and B.~Kim, 
and I.~Ciocan--Fontanine, that
$$QH^*(\Fl )\simeq\overline P_n/\wt I_n, 
$$
and the residue pairing defined by $\wt I_n$ may be naturally identified 
with the intersection form on the quantum cohomology ring. We will call 
this residue pairing as quantum pairing on $\overline P_n/\wt I_n$ and denote 
it by $\langle ,\rangle_Q$.

\section{Quantum double Schubert polynomials.}

Quantum double Schubert polynomials are closely related with the 
equivariant quantum cohomology. Let us remind the result of A.~Givental 
and B.~Kim [GK] (see also [K2]) on the structure of the equivariant 
quantum cohomology algebra of the flag variety $Fl_n$:
$$QH^*_{U_n}(Fl_n)\cong{\bf Z}[x_1,\ldots ,x_n,y_1,\ldots ,y_n,q_1,\ldots 
,q_{n-1}]/\wt J,
$$
where the ideal $\wt J$ generated by 
$$e_i(x_1,\ldots ,x_n~|~q_1,\ldots ,q_{n-1})-e_i(y_1,\ldots ,y_n), \ 1\le 
i\le n.
$$

In classical case $q=0$, the double Schubert polynomials $\s_w(x,y)$
represent the equivariant cohomology classes [F]. Quantum double Schubert
polynomials have to play the similar role for the quantum equivariant
cohomology ring. Let us define at first the ``top'' quantum double
Schubert polynomial $\wt{\s}_{w_0}(x,y)$.

Let $x=(x_1,\ldots ,x_n)$, $y=(y_1,\ldots ,y_n)$ be two sets of 
variables, put
$$\wt{\s}_{w_0}(x,y):=\wt{\s}_{w_0}^{(q)}(x,y)=
\prod_{i=1}^{n-1}\Delta_i(y_{n-i}~|~x_1,\ldots ,x_i),
$$
where $\Delta_k(t~|~x_1,\ldots ,x_k):=\ds\sum_{j=0}^kt^{k-j}
e_j(x_1,\ldots ,x_k~|~q_1,\ldots ,q_{k-1})$ is the generating function 
for the quantum elementary symmetric polynomials in $x_1,\ldots ,x_k$,
i.e. ~~
$ \Delta_k(t|x):=\ds\sum_{i=1}^ke_i(x|q)t^i=$
\begin{equation}
\det\left(\begin{array}{ccccccc} x_1+t & q_1&0 &\ldots &\ldots 
&\ldots &0\\
-1 & x_2+t & q_2 &0 &\ldots &\ldots & 0\\
0 & -1 &x_3+t & q_3 & 0 &\ldots & 0 \\
\vdots &\ddots &\ddots & \ddots &\ddots &\ddots &\vdots \\
0&\ldots & 0 &-1&x_{k-2}+t &q_{k-2} & 0 \\
0 &\ldots &\ldots & 0 &-1 & x_{k-1}+t & q_{k-1}\\
0 &  \ldots &\ldots &\ldots & 0 & -1 & x_k+t
\end{array}\right) .
\end{equation}
\vskip 0.5cm

\begin{de} For each permutation $w\in S_n$, the quantum double Schubert 
polynomial $\wt{\s}_w(x,y)$ is defined to be
$$\wt{\s}_w(x,y)=\partial_{ww_0}^{(y)}\wt{\s}_{w_0}(x,y),
$$
where divided difference operator $\partial_{ww_0}^{(y)}$ acts on the $y$ 
variables.
\end{de}
{\bf Remark 2} \  $i)$ In the "classical limit" $q_1=\cdots =q_{n-1}=0$,
$$\wt{\s}_w(x,y)|_{q=0}=\partial_{ww_0}^{(y)}\s_{w_0}(x,y)=\s_{w^{-1}}(y,x)=
\s_w(x,y),
$$
i.e. $\wt{\s}_w(x,y)|_{q=0}=\s_w(x,y)$.

$ii)$ (Stability) Let $m>n$ and let $i:S_n\hookrightarrow S_m$ be the 
embedding. Then
$$\wt{\s}_w(x,y)=\wt{\s}_{i(w)}(x,y).
$$

$iii)$ One can check that the ring
$${\bf Z}[x_1,\ldots ,x_n,y_1,\ldots ,y_n,q_1,\ldots ,q_{n-1}]/{\wt J}
$$
is a free module of dimension $n!$ over the quotient ring $\wt R$ with basis 
either $\wt{\s}_w(x)$ or $\wt{\s}_w(x,y)$, $w\in S_n$, where
$$\wt R:=\frac{{\bf Z}[x_1,\ldots ,x_n,q_1,\ldots ,q_{n-1}]\otimes 
{\bf Sym}[y_1,\ldots ,y_n]}{\wt J}.
$$
\smallskip \\ 
{\bf Example.} Quantum double Schubert polynomials for $S_3$:
\begin{eqnarray*}\wt{\s}_{s_1s_2s_1}(x,y) &=& 
(x_1+y_2)(x_1+y_1)(x_2+y_1)+q_1(x_1+y_2),\\
\wt{\s}_{s_2s_1}(x,y) &=& (x_1+y_1)(x_1+y_2)-q_1,\\
\wt{\s}_{s_1s_2}(x,y) &=& (x_1+y_1)(x_2+y_1)+q_1,\\
\wt{\s}_{s_1}(x,y) &=& x_1+y_1,\\
\wt{\s}_{s_2}(x,y) &=& x_1 +x_2+y_1+y_2,\\
\wt{\s}_{\rm id}(x,y) &=& 1.
\end{eqnarray*}

For the list of the quantum double Schubert polynomials corresponding
to the symmetric group $S_4$, see appendix A.

{\begin{th} Let $z=(z_1,\ldots ,z_n)$ be a third set of variables. Then
\begin{equation}
\langle\wt{\s}_{w_0}(x,y),\wt{\s}_{w_0}(x,z)\rangle_Q^{(x)}=C(y,z),\label{6}
\end{equation}
where the upper index $x$ means that the quantum pairing is taken in the 
$x$ variables, and
$$C(x,y)=\sum_{w\in S_n}\s_w(x)\s_{w_0w}(y)
$$
is the "canonical" element in the tensor product $H^*(\Fl )\otimes 
H^*(\Fl )$.
\end{th}

Theorem 1 plays the important role in our approach to the 
quantum Schubert polynomials. We will give a proof later, and now let us 
consider some applications of the formula (10). 

\section{Quantum Schubert polynomials.}

\subsection{Definition.}

Let us remind the result of A.~Givental and B.~Kim, and 
I.~Ciocan-Fontanine on the structure of the small quantum cohomology ring 
of flag variety $\Fl$ 
$$QH^*({\Fl})\cong{\bf Z} [ x_1,\ldots ,x_n,q_1,\ldots ,q_{n-1} ] /\wt I,
$$
where the ideal $\wt I$ is 
generated by the quantum elementary symmetric polynomials $\wt e_i(x):=
e_i(x_1,\ldots ,x_n|q_1,\ldots ,q_{n-1})$, $1\le i\le n$ with generating 
function $\Delta_n(t|x)$, see (9).

We define a pairing on the ring of polynomials ${\bf Z}[x;q]$ 
and the quantum cohomology
ring $QH^*(\Fl )\simeq{\bf Z}[x;q]/\wt I$ using 
the Grothendieck residue
$$\langle f,g\rangle_Q=\Res_{\wt I}(fg),~f,g\in{\bf Z}[x_1,\ldots 
,x_n,q_1,\ldots ,q_{n-1}].
$$
Then

1) $\langle f,g\rangle_Q=0$ if $f\in \wt I$; 

\vskip 0.2cm

2) $\langle f,g\rangle_Q$ defines a nondegenerate pairing in $QH^*(\Fl )$.

\begin{de}\label{102}~
Define the quantum Schubert polynomials $\wt{\s}_w:=\wt{\s}_w(x)$ 
as Gram--Schmidt's 
orthogonalization  of the set of lexicographically ordered monomials
\hbox{$\{x^I~|~I\subset\delta\}$} with respect to the quantum residue pairing
$\langle f,g\rangle_Q$:

1) $\langle\wt{\s}_u,\wt{\s}_v\rangle_Q=\langle\s_u,\s_v\rangle 
=\left\{\begin{array}{ll} 1,& \mbox{\rm{if} $v=w_0 u$}\\ 0,
& \mbox{\rm{otherwise}}\end{array}\right.$

2) $\wt{\s}_w(x)=x^{c(w)}+\ds\sum_{I<c(w)}a_I(q)x^I$, where 
$a_I(q)\in{\bf Z}[q_1,\ldots ,q_{n-1}]$ and\break $I<c(w)$ means the 
lexicographic order.
\end{de}
Here $c(w)$ is the code of a permutation $w\in S_n$, [M1], p.9.\medskip \\
{\bf Remark 3} \ This definition is the analogue of the 
characterization of Schubert polynomials from Proposition 4. 
\smallskip \\ 
{\bf Example.} For the symmetric group $S_3$, we have
$$\langle x_1^2x_2,x_1^2\rangle_Q=q_1,~~\langle 
x_1^2x_2,x_1x_2\rangle_Q=-2q_1.
$$
Consequently,
$$\wt{\s}_1=x_1,~\wt{\s}_2=x_1+x_2,~\wt{\s}_{12}=x_1x_2+q_1,~\wt{\s}_{21}=
x_1^2-q_1,~\wt{\s}_{121}=x_1^2x_2+q_1x_1.
$$
Let us remark that in our example ($n=3$) 
$\wt{\s}_{121}=\wt{\s}_{w_0^{(3)}}(x)=\wt e_1(x_1)\wt e_2(x_1,x_2)$. More 
generally, we have

\begin{pr} Let $w_0\in S_n$ be the longest element. Then
$$\wt{\s}_{w_0}(x)=\wt e_1(x_1)\wt e_2(x_1,x_2)\ldots \wt 
e_{n-1}(x_1,\ldots x_{n-1}).
$$
\end{pr}

In other words, the quantum Schubert polynomial corresponding to the 
longest element of the symmetric group $S_n$, is equal to the product of 
all principal minors of the Jacobi matrix $\ds\left(\frac{\partial
e_i(x|q)}{\partial x_j}\right)_{1\le i,j\le n}$. We can also compute the 
Grothendieck residue w.r.t. ideal $\wt I$ of the Jacobian 
$\det\ds\left(\frac{\partial e_i(x|q)}{\partial x_j}\right)_{1\le 
i,j\le n}$.

\begin{pr} (cf. [EL])\\
$$\det\left(\frac{\partial e_i(x|q)}{\partial x_j}\right)\equiv 
n!\wt{\s}_{w_0}(x)\ (\mod \wt I),
$$
where $e_i(x|q)=e_i(x_1,\ldots ,x_n~|~q_1,\ldots ,q_{n-1})$, $1\le i\le 
n$, are the quantum elementary symmetric functions.
\end{pr}

Remind that $n!$ is equal to the Euler number of $\Fl $.

\subsection{Orthogonality.}

We use the Jack--Macdonald type definition ([M2], Chapter VI) of the
quantum Schubert polynomials, see Definition 5. On this way the 
orthogonality of quantum Schubert polynomials is valid by ``definition''.
We are going to prove that the $y=0$ specialization of quantum double
Schubert polynomials $\wt{\s}_w(x,0)$ also satisfies the conditions 1) and
2) of Definition 5. As a corollary, we obtain that the specialization 
$\wt{\s}_w(x,0)$ coincides with the quantum Schubert polynomial
$\wt{\s}_w(x)$ from Definition 5.
\begin{th}\label{101} Let $v,w\in S_n$. Then
$$\langle \wt{\s}_v(x,0),\wt{\s}_w(x,0)\rangle_Q=\left\{\begin{array}{ll}
1,& \mbox{{\rm if} $w=w_0v$},\\ 0,& \mbox{\rm otherwise}.\end{array}\right.
$$
\end{th}

{\it Proof.} Let us apply the operator 
$\partial_{vw_0}^{(y)}\partial_{ww_0}^{(z)}$ to the both sides of 
(\ref{6}). The LHS gives
$$\partial_{vw_0}^{(y)}\partial_{ww_0}^{(z)}\langle\wt{\s}_{w_0}(x,y),
\wt{\s}_{w_0}(x,z)\rangle_Q^{(x)}=\langle\wt{\s}_v(x,y),
\wt{\s}_w(x,z)\rangle_Q^{(x)}.
$$
The RHS transforms to the following form $\ds\sum_{u\in 
S_n}\partial_{vw_0}^{(y)}\s_u(y)\partial_{ww_0}^{(z)}\s_{w_0u}(z)$. Now 
taking $y=z=0$ we obtain an equality
\begin{equation} \langle\wt{\s}_v(x,0)\wt{\s}_w(x,0)\rangle_Q=\sum_{u\in 
S_n}\eta (\partial_{vw_0}\s_u)\eta (\partial_{ww_0}\s_{w_0u}),\label{7}
\end{equation}
where $\eta :P_n\to{\bf Z}$ is the homomorphism defined by $\eta 
(x_i)=0$ $(1\le i\le n)$. It is clear that
$$\eta (\partial_v\s_u)=\left\{\begin{array}{ll} 1, &{\rm if} \ v=u,\\
0, & {\rm otherwise}. \end{array}\right.
$$
Thus, the RHS of (\ref{7}) is equal to 1 if $w_0ww_0=vw_0$ and is equal 
to 0 otherwise.

\qed\smallskip \\
{\bf Remark 4} \ Orthogonality of quantum Schubert 
polynomials was proven in [FGP], using a combinatorial definition, see
{\it ibid}, Section 5; the proof is highly non trivial.

\subsection{Quantum Cauchy formula.}

\begin{th} Let $\wt{\s}_w(x):=\wt{\s}_w(x,0)$, then
\begin{equation}\sum_{w\in S_n}\wt{\s}_w(x)\s_{ww_0}(y)=\wt{\s}_{w_0}(x,y).
\label{8}
\end{equation}
\end{th}

{\it Proof.} Let us apply the divided difference operator 
$\partial_{ww_0}^{(z)}$ to the both sides of (\ref{6}) and then
take $z=0$. The right hand side transforms to the following form
$$\sum_{u\in 
S_n}\s_u(y)\partial_{ww_0}^{(z)}\s_{w_0u}(z)|_{z=0}=\s_{w_0ww_0}(y).
$$
As for the LHS, it takes the form 
$\langle\wt{\s}_{w_0}(x,y),\wt{\s}_w(x)\rangle_Q$. Hence,
$$\langle\wt{\s}_{w_0}(x,y),\wt{\s}_w(x)\rangle_Q=\s_{w_0ww_0}(y).
$$
The last identity is equivalent to (\ref{8}).
\qed

More generally, we have
\begin{pr}\label{104}
\begin{eqnarray}\sum_{w\in S_n}\wt{\s}_w(x,z)\s_{ww_0}(y,-z) &=&
\wt{\s}_{w_0}(x,y),\\
\sum_{u\in S_n,\ l(u)+l(uw^{-1})=l(w)}\wt{\s}_u(x,z)\s_{uw^{-1}}(y,-z) & = &
\wt{\s}_w(x,y).
\end{eqnarray}
\end{pr}

{\it Proof.} Let us apply the Interpolation formula to 
$f(x)=\wt{\s}_{w_0}(x,y)$ and then divided difference operator 
$\partial_{ww_0}^{(y)}$.

\qed
\begin{col}
$$C^{(q,q')}(x,y):=\sum_{w\in S_n}\wt{\s}_w^{(q)}(x)
\wt{\s}_{w_0w}^{(q')}(y)=\langle\wt{\s}_{w_0}^{(q)}(x,z),
\wt{\s}_{w_0}^{(q')}(y,z)\rangle^{(z)},
$$
where the upper index $z$ means that the scalar product is taken in the
$z$ variables.
\end{col}

\begin{col}\label{107}
$$\sum_{w\in S_n}\Omega_w\Omega_w^*=C^{(q,q)}(x,x).
$$
\end{col}
One can show that
$$C^{(q,q)}(x,x)=\sum_{w\in S_n}\wt{\s}_w(x)\wt{\s}_{ww_0}(x)\equiv n!\:
\wt{\s}_{w_0}(x)~(\mod\wt I).
$$

Let us summarize our results. It follows from Theorem \ref{101} 
that polynomials 
$\wt{\s}_w(x,0)$ are orthogonal with respect to the quantum pairing $\langle 
,\rangle_Q$. It is also clear that $\wt{\s}_w(x,0)|_{q=0}=\s_w(x)$ and
$\wt{\s}_w(x)=x^{c(w)}+$lower degree terms w.r.t. lexicographic order 
on the set of monomials. These 
two properties characterize the polynomials $\wt{\s}_w(x,0)$ uniquely, 
consequently, the polynomials $\wt{\s}_w(x,0)$ coincide with the quantum 
Schubert polynomials $\wt{\s}_w(x)$ from Definition \ref{102}. As a matter of 
fact, we obtain the Lascoux-Sch\"utzenberger type formula for quantum 
Schubert polynomials.

\begin{th}\label{t7} Let $w\in S_n$, then
$$\wt{\s}_w(x)=\partial_{ww_0}^{(y)}\wt{\s}_{w_0}(x,y)|_{y=0}.
$$
\end{th}

\section{Quantization.}

\subsection{Definition.}

Let $f\in P_n={\bf Z}[x_1,\ldots ,x_n]$ be a polynomial. According to the
Interpolation formula,
$$f(x)=\sum_{w\in S^{(n)}}\partial_w^{(y)}f(y)\s_w(x,y).
$$
We define a quantization $\wt f$ of the function $f$ by the rule
\begin{equation}\wt f(x)=\sum_{w\in S^{(n)}}\partial_w^{(y)}f(y)
\wt{\s}_w(x,y)|_{\overline P_n},\label{51}
\end{equation}
where for a polynomial $f\in\overline P_{\infty}$, the symbol 
$f\mid_{\overline P_m}$ means the restriction of $f$ to the ring of
polynomials $\overline P_m$,
i.e. the specialization $x_{m+1}=x_{m+2}=\cdots =0$ and $q_m=q_{m+1}=
\cdots =0$.

Hence, the quantization is a ${\bf Z}[q_1,\ldots ,q_{n-1}]$--linear 
map $P_n\to\overline P_n$.

The main property of quantization is that it preserves the pairings, i.e.
\begin{equation}\langle\wt f,\wt g\rangle_Q=\langle f,g\rangle , \ 
f,g\in P_n. \label{110}
\end{equation}
It follows from (\ref{110}) that the quantization map maps the ideal 
$I_n\subset\overline P_n$ into ideal $\wt I_n\subset\overline P_n$.
\smallskip \\ 
{\bf Remark 5} \ $i)$ Quantization does not preserve multiplication, i.e. in 
general $\wt f\cdot \wt g\ne \wt{fg}$. For example, if 
$f=\ds\sum_{i=1}^n\alpha_ix_i$ is a linear form, then
(quantum Monk's formula, see [FGP] and our Section 9)
$$\wt f\wt{\s}_w-\wt{f\s_w}=\sum (\ld_i-\ld_j)q_{ij}\wt{\s}_{wt_{ij}},
$$
summed over $i<j$ such that $l(w)=l(wt_{ij})+l(t_{ij})$. Here 
$q_{ij}=q_iq_{i+1}\ldots q_{j-1}$.

$ii)$ It is clear that if $f\in H_n$, then $\wt f\in\overline H_n$.
 
$iii)$ It follows from Proposition \ref{104}, that the quantum double Schubert 
polynomials $\wt{\s}_w(x,y)$ are the quantization of classical ones.

$iv)$ It follows from Interpolation formula and quantization procedure, that

$\bullet$ Quantum Schubert polynomials $\wt{\s}_w(x)$, $w\in S_n$ form a 
$\overline I$-basis in $\overline P_n$.

$\bullet$ Quantum Schubert polynomials $\wt{\s}_w(x)$, $w\in S^{(n)}$ form a 
${\bf Z}[q_1,\ldots ,q_{n-1}]$--basis of $\overline P_n$.

$\bullet$ Quantum Schubert polynomials $\wt{\s}_w(x)$, $w\in S_n$ form a 
${\bf Z}[q_1,\ldots ,q_{n-1}]$--basis of $\overline H_n=
H_n\otimes{\bf Z}[q_1,\ldots ,q_{n-1}]$.

The proof of the statement $iv)$ can be found in [FGP].
\smallskip

Now we are going to describe another families of polynomials having a nice 
quantization.

\subsection{Elementary and complete polynomials.}

Let $\delta :=\delta_n=(n-1,n-2,\ldots ,1,0)$, and consider the set $\t $ 
of sequences $I=(i_1,\ldots ,i_n)\in{\bf Z}^n$ such that $0\le i_j\le 
n-j$ for all $j=1,\ldots ,n$. It is clear that $|\t |=n!$, and

$\bullet$ $P_n$ is a free $\Lambda_n$-module of rank $n!$ with basis 
$\{x^I=x_1^{i_1}x_2^{i_2}\ldots x_n^{i_n}~|~I\in\t\}$. 

Follow to [LS2], for each $I\in\t$ let us define the elementary polynomial 
$e_I(x)$ as the following product
$$\prod_{k=1}^{n-1}e_{i_k}(x_1,\ldots ,x_{n-k}).
$$

$\bullet$ (Lascoux-Sch\"utzenberger [LS2]) $P_n$ is a free 
$\Lambda_n$-module of rank $n!$ with basis $\{e_I(x)~|~I\in\t\}$.

\begin{de} For each sequence $I\in\t$ 
the quantum elementary polynomial ${\wt e}_I(x)$ is defined to be
$${\wt e}_I(x)=\prod_{k=1}^{n-1}{\wt e}_{i_k}(x_1,\ldots ,x_{n-k}),
$$
where ${\wt e}_k(x_1,\ldots ,x_m):=e_k(x_1,\ldots ,x_m~|~q_1,\ldots 
,q_{m-1})$ are the quantum elementary symmetric functions.
\end{de}

\begin{th}\label{t8} Assume that $I\subset\delta$. Then ${\wt e}_I(x)$ is the 
quantization of elementary polynomial $e_I(x)$.
\end{th}

{\it Proof.} It is enough to prove the following 
\begin{pr}\label{10} If  $I\subset\delta$, then
\begin{equation} {\wt e}_I(x)=\sum_{w\in S_n}\wt{\s }_w(x)\eta 
(\partial_we_I).\label{9}
\end{equation}
\end{pr}
We are going to show that Proposition \ref{10} follows from the quantum 
Cauchy formula. First of all, let us remark that
\begin{eqnarray*}\wt{\s}_{w_0}(x,y)&=&\sum_{I\subset\delta}{\wt 
e}_I(x)y^{\delta -I},\\
\wt{\s}_{w_0}(x,z)&=&\sum_{w\in S_n}\wt{\s}_{w_0ww_0}(x)\s_{w_0w}(z).
\end{eqnarray*}
Substituting these two expressions in (\ref{9}), we obtain a formula for 
the classical Schubert polynomials
\begin{equation}\s_w(y)=\sum_{I\subset\delta}\langle{\wt 
e}_I(x),\wt{\s}_{w_0ww_0}(x)\rangle_Qy^{\delta -I}.\label{11}
\end{equation}
It follows from (\ref{11}) that
\begin{equation}\langle{\wt e}_I(x),\wt{\s}_{w_0ww_0}(x)\rangle_Q=
\langle e_I(x),\s_{w_0ww_0}(x)\rangle .\label{12}
\end{equation}
Conversely, the quantum Cauchy formula (\ref{8}) follows from the 
classical one and (\ref{12}). To continue, let us remark that
\begin{eqnarray*}
\langle e_I(x),\s_{w_0ww_0}(x)\rangle &=&
\langle e_I(x),\partial_{w_0w^{-1}}\s_{w_0}(x)\rangle =\\
\langle\partial_{ww_0}e_I(x),\s_{w_0}(x)\rangle 
&=&\eta (\partial_{ww_0}e_I(x)).
\end{eqnarray*}
As a corollary we obtain a formula for Schubert polynomials, 
which seems to be new,
\begin{equation} \s_w(x)=\sum_{I\subset\delta}\eta 
(\partial_{ww_0}e_I(x))x^{\delta -I}.\label{13}
\end{equation}
Formula (\ref{13}) gives a geometric interpretation of the coefficients 
$a_{I,w}$ of the Schubert polynomial 
$\s_w(x)=\sum_{I\subset\delta}a_{I,w}x^{\delta -I}$ as the intersection 
numbers
$$a_{I,w}=\eta (\partial_{ww_0}e_I(x))=\langle 
e_I(x),\s_{w_0ww_0}(x)\rangle \ge 0.
$$
Finally, let us finish a proof of Theorem \ref{t8}. We have
\begin{eqnarray*}\sum_{I\subset\delta}{\rm RHS}(\ref{9})y^{\delta -I}&=&
\sum_{w\in S_n}\wt{\s}_w(x)\sum_{I\subset\delta}\eta 
(\partial_we_I)y^{\delta -I}\\
=\sum_{w\in 
S_n}\wt{\s}_w(x)\s_{ww_0}(y)&=&\wt{\s}_{w_0}(x,y)=\sum_{I\subset\delta}{\wt 
e}_I(x)y^{\delta -I}.
\end{eqnarray*}
Hence, ${\rm RHS}(\ref{9})={\wt e}_I(x)$.
\qed

\begin{col} If $I$ and $J$ belong to $\t$, then
$$\langle{\wt e}_I(x),{\wt e}_J(x)\rangle_Q=\langle 
e_I(x),e_J(x)\rangle .
$$
\end{col}
\vskip 0.3cm
{\bf Remark 6} \ In the next section we will give a proof of Corollary 3 
using a geometric technique due to I.~Ciocan--Fontanine [C] (see also 
[K1]). Repeating our arguments in the reverse order, we see that the 
quantum Cauchy formula (12), as well as Theorems 1 and 3, follow directly from 
Corollary 3.
\medskip

Using quantum Cauchy formula (\ref{8}), we can describe a transition matrix 
between quantum Schubert polynomials and quantum elementary polynomials.

\begin{th} 
$$\wt{\s}_w(x)=\ds\sum_{I\subset\delta}\wt e_I(x)\eta 
(\partial_{ww_0}x^{\delta -I}).$$
\end{th}

{\it Proof.} It follows from Cauchy's formula that
$$\sum_{w\in S_n}\wt{\s}_w(x)\s_{ww_0}(y)=
\sum_{I\subset\delta}\wt e_I(x)y^{\delta -I}.
$$
Consequently,
$$\wt{\s}_{w_0ww_0}(x)=\sum_{I\subset\delta}\wt e_I(x)\langle y^{\delta 
-I},\s_w(y)\rangle .
$$
Now we have
$$\langle y^{\delta -I},\s_w(y)\rangle =\langle y^{\delta -I},
\partial_{w^{-1}w_0}\s_{w_0}(y)\rangle =\langle\partial_{w_0w}y^{\delta 
-I},\s_{w_0}(y)\rangle =\eta (\partial_{w_0w}y^{\delta -I}).
$$
\qed
\smallskip \\ 
{\bf Example.} Take the permutation $w=24531\in S_5$. It is easy to check 
that 
$ww_0=42135=s_1s_2s_1s_3$ and there exists 6 monomials 
$x^I$ such that $I\subset (43210)$ and $\eta (\partial_{1213}x^I)\ne 0$. 
They are
$$\begin{array}{rcccccc} x^I:& x_1^2x_2x_3,& x_1^2x_2x_4, & x_1^2x_3^2,&
x_1x_2^2x_3,& x_1x_2^2x_4,& x_2^2x_3^2\\ \\
\eta(\partial_{1213}x^I):& +1& -1& -1& -1& +1& +1
\end{array}
$$

We can check using Theorem \ref{t7}, that
$$\wt{\s}_{24531}(x)=\wt e_{2211}(x)-\wt e_{2220}(x)-\wt e_{2301}(x)-
\wt e_{3111}(x)+\wt e_{3120}(x)+\wt e_{4101}(x).
$$

Now let us consider a problem how to quantize monomials.

\begin{pr} Let $I\in\t$, then
\begin{eqnarray*} x^I&=&\sum_{w\in S_n, \  l(w)=|I|}\eta 
(\partial_wx^I)\s_w(x),\\ \wt x^I&=&\sum_{w\in S_n, \ l(w)=|I|}\eta 
(\partial_wx^I)\wt{\s}_w(x).
\end{eqnarray*}
\end{pr}
\begin{col} Let $v,w\in S_n$, then
$$\sum_{I\subset\delta}\eta (\partial_vx^I)\eta (\partial_{ww_0}e_{\delta 
-I}(x))=\delta_{v,w}.
$$
\end{col}
\begin{col}
$$\wt{\s}_w(x)=\sum_{I\subset\delta}\eta (\partial_{ww_0}e_{\delta 
-I}(x))\wt x^I.
$$
\end{col}
\begin{col}
$$\wt{\s}_{w_0}(x,y)=\sum_{I\subset\delta}\wt x^Ie_{\delta -I}(y).
$$
\end{col}

Now let us consider a problem how to quantize the complete homogeneous 
symmetric functions 
\[ h_k^m := \sum_{i_1+\cdots+ i_m =k}x_1^{i_1}\cdots x_m^{i_m} . \] 

We define the quantum complete homogeneous symmetric function 
$\wt{h}_k^m = \wt{h}_k(x_1,\ldots ,x_m)$  of degree $k$ using the 
generating function 
$$H_m(t\mid x) := H_m(t\mid x_1,\ldots x_m) = \left( t^m \Delta_m(t^{-1}
\mid -x) \right) ^{-1}=\sum_{k\geq 0}t^k \wt{h}_k^m. 
$$
It is clear that 
\[ t^m \Delta_m(t^{-1}\mid x) = \sum_{k\geq 0} t^k \wt{e}_k^m , \] 
and 
\[ \Delta_m(t^{-1}\mid x)\cdot H_m(t\mid -x) =1. \] 
The last relation gives possibility to express the quantum complete 
homogeneous symmetric functions in terms of quantum elementary ones: 
\begin{eqnarray}
\wt{h}_k^m & = & {\rm det} 
\left( \wt{e}_{j-i+1}^m \right)_{1\leq i,j\leq k} \nonumber \\ 
 & = & \left({\rm det} 
\left( \wt{e}_{j-i+1}^{m+k-i} \right)_{1\leq i,j\leq k}
\right)\mid_{_{\overline P_m}}.
\end{eqnarray} 

The equality (20) follows from the recurrence relation for 
$\wt{e}_k^N$ (see [C], [GK]): 
\[ \wt{e}_k^N = \wt{e}_k^{N-1} + x_N \wt{e}_{k-1}^{N-1} 
+ q_{N-1} \wt{e}_{k-2}^{N-2}. \] 
But each term in the expansion of the determinant (20) is a quantum elementary 
polynomial in $\overline P_{m+k}$. Hence, it follows from Theorem 7 that 
the quantum complete 
homogeneous symmetric function $\wt{h}_k^m$ is the quantization of 
classical one. 

Finally, let us define the complete and quantum complete polynomials.
\begin{de} For each sequence $I\subset\delta_n$ the complete and quantum 
complete polynomials $h_I(x)$ and $\wt h_I(x)$ are defined to be
\begin{eqnarray} h_I(x) & = & \prod_{k=1}^{n-1}h_{i_k}(x_1,\ldots ,x_k),\\
\wt h_I(x) & = & \prod_{k=1}^{n-1}\wt h_{i_k}(x_1,\ldots ,x_k).
\end{eqnarray}
\end{de}
{\bf Remark 8} \ It follows from (20) that if $m+k>n$ then 
$\wt h_k^m\in\wt I_n$.
Hence, if $I\not\subset\delta$, then $h_I(x)\in I_n$ and $\wt h_I(x)\in
\wt I_n$.

It is not difficult to see that $P_n$ is a free $\Lambda_n$-module of rank
$n!$ with basis $\{ h_I(x)~|~I\in \t\}$.\smallskip \\
{\bf Theorem 5'} {\it If $I\subset\delta$, then $\wt h_I(x)$ is the 
quantization of complete polynomial $h_I(x)$.}

{\it Proof.} See Remark 10 in Section 7. \qed

\subsection{Canonical involution $\omega$.}

There exists an involution $\omega$ of the ring $\overline P_n[y]$ given
by $\omega (x)={\covec x}$, $\omega (y)={\covec y}$, $\omega (q)=
{\covec q}$, where for any sequences $z=(z_1,\ldots , z_m)$ we define 
$\covec z$ to be equal to $(z_m, z_{m-1},\ldots ,z_1):=\covec z$. It 
is clear from the definition of quantum elementary symmetric functions
 $e_i(x|q)$, see (9), Section 3, that
$$\omega (e_i(x|q))=e_i(x|q)
$$
and thus the involution $\omega$ preserves the ideal $\wt I_n$ (as well as
the ideals $I_n$, $J_n$ and $\wt J_n$).
\begin{pr} 
$$\omega (\wt{\s}_u^{(q)}(x,y))\equiv\epsilon (u)
\wt{\s}_{w_0uw_0}^{(q)}(x,y)~\mod\wt J_n,
$$ 
where $\epsilon (u)=(-1)^{l(u)}.$
\end{pr}

{\it Proof.} First of all, if $u=w_0$, then
\begin{equation}
\omega (\wt{\s}_{w_0}^{(q)}(x,y))\equiv\epsilon (w_0)
\wt{\s}_{w_0}^{(\covec q)}(\covec x,\covec y)~\mod \wt J_n.\label{21}
\end{equation}
But $\omega\partial_u=\epsilon (u)\partial_{w_0uw_0}\omega$ (see [M1], (2.12)).
Thus, applying the divided difference operator $\epsilon (u)
\partial_{w_0uw_0}^{(y)}$
to the both sides of (\ref{21}), we obtain
$$\omega (\wt{\s}_{uw_0}^{(q)}(x,y))\equiv\epsilon (w_0)
\wt{\s}_{w_0u}^{(\covec q)}(\covec x,\covec y)~\mod \wt J_n.
$$
\qed

Finally, let us describe the action of involution $\omega $ on the 
elementary polynomials.
\begin{pr} 
$$\omega (\wt e_I(x))\equiv (-1)^{|I|}
\wt h_{\covec I}(x)~\mod \wt I_n.
$$
\end{pr}\bigskip 
{\bf Remark 9} \ To our knowledge, originaly, construction of the 
quantization map, using a remarkable family of commuting operators $X_i$,
appeared in [FGP]. We use a different definition of quantization map, but
it can be shown that two forms of quantizations are equivalent.
For original proofs of Theorem 5 and 5', and Proposition 15,
see Corollary 4.6, Corollary 7.16 and Proposition 7.13 in [FGP]. 

\section{Quantum cohomology ring of flag variety.}

Quantum cohomology ring of the flag variety ${\Fl}$ is a
deformed ring of 
the ordinary cohomology ring $H^*({\Fl} ,{\bf Z})$. The structure constants of 
the quantum cohomology ring are given by the Gromov--Witten invariants. Let 
$\Omega_{w_1},\ldots ,\Omega_{w_m}$ ($w_i\in S_n$) be Schubert cycles. We 
denote by $M_{\bar d}({\bf P}^1,Fl_n)$ the moduli space of morphisms from 
${\bf P}^1$ to $\Fl$ of multidegree $\bar{d}=(d_1,\ldots ,d_{n-1})$. We 
consider the restriction of the universal map for $t\in{\bf P}^1$:
$$ev_t~:~M_{\bar d}({\bf P}^1,\Fl )\times\{ t\}\hookrightarrow 
M_{\bar d}({\bf P}^1,\Fl )\times{\bf P}^1{\buildrel
ev\over\rightarrow}
\Fl ,~~(f,p)\mapsto f(p).
$$
Let $\Omega_w(t)=ev_t^{-1}(\Omega_w)$.

\begin{th} \mbox{(I.~Ciocan-Fontanine)}. If
$\ds\sum_{i=1}^ml(w_i)=\frac{n(n-1)}{2}
+2\sum d_i$ and\break $t_1,\ldots ,t_m\in{\bf P}^1$ are distinct, then for 
general translates of ~$\Omega_{w_i}$, the number of points in 
$\ds\bigcap_{i=1}^m\Omega_{w_i}(t_i)$ is finite and independent of $t_i$ and 
the translates of $\Omega_{w_i}$.
\end{th}

\begin{de} The Gromov--Witten invariant is defined as an intersection 
number
$$\langle\Omega_{w_1}\ldots\Omega_{w_m}\rangle_{\bar d}=\left\{
\begin{array}{ll}\#\ds\bigcap_i\Omega_{w_i}(t_i), & \mbox{{\rm if} 
$\sum l(w_i)=\frac{n(n-1)}{2}+2\sum d_i$}\\ 0, & \mbox{\rm otherwise}
\end{array}\right. .
$$
\end{de}

Now we can define the quantum multiplication as a linear map
$$m_q~:~{\rm Sym~}(H^*(\Fl ,{\bf Z})[q_1,\ldots ,q_{n-1}])\to H^*(\Fl ,{\bf 
Z})[q_1,\ldots ,q_{n-1}]
$$
given by
$$m_q(\prod_{i=1}^m\Omega_{w_i})=\sum_{\bar d}q^{\bar d}
\sum_w\langle\Omega_w\Omega_{w_1}\cdots\Omega_{w_m}\rangle_{\bar d}
\Omega^*_w,
$$
where $q^{\bar d}=q^{d_1}_1\cdots q_{n-1}^{d_{n-1}}$ and $(\Omega^*_w)$ is 
the dual basis of $(\Omega_w)$.

Then the quantum cohomology ring $QH^*(\Fl )$ is a commutative and 
associative ${\bf Z}[q_1,\ldots ,q_{n-1}]$ -- algebra.

Let $0=E_0\subset E_1\subset\cdots\subset E_n={\bf C}^n\otimes{\cal O}_F$
be the universal flag of subbundles on $\Fl$.

\begin{th} \label{106}(A.~Givental and B.~Kim, I.~Ciocan--Fontanine).\\
The small quantum cohomology 
ring is generated by $x_i=c_1(E_{n-i+1}/E_{n-i})$, $i=1,\ldots ,n$, as 
a ${\bf Z}[q_1,\ldots ,q_{n-1}]$-algebra and
$$QH^*(\Fl )\cong{\bf Z}[x_1,\ldots ,x_n,q_1,\ldots 
,q_{n-1}]/(e_1(x|q),\ldots ,e_n(x|q)),
$$
where $e_i(x|q)$ is given by the expansion of the following determinant
$$\det\left(\begin{array}{ccccccc} x_1+t & q_1&0 &\ldots &\ldots 
&\ldots &0\\
-1 & x_2+t & q_2 &0 &\ldots &\ldots & 0\\
0 & -1 &x_3+t & q_3 & 0 &\ldots & 0 \\
\vdots &\ddots &\ddots & \ddots &\ddots &\ddots &\vdots \\
0&\ldots & 0 &-1&x_{n-2}+t &q_{n-2} & 0 \\
0 &\ldots &\ldots & 0 &-1 & x_{n-1}+t & q_{n-1}\\
0 &  \ldots &\ldots &\ldots & 0 & -1 & x_n+t
\end{array}\right) 
$$
$$=t^n+e_1(x|q)t^{n-1}+\cdots +e_n(x|q).
$$
\end{th}

It follows from Theorem \ref{106} that any Schubert cycle $\Omega_w$ 
may be expressed as a 
polynomial $\wh{\s}_w(x ,q)$ in $QH^*(\Fl )$. The polynomial 
$\wh{\s}_w(x ,q)$ is a deformation of the Schubert polynomial 
$\s_w(x )$ and $\wh{\s}_w(x ,0)=\s_w(x)$. Consider the correlation 
function
$$\langle\Omega_{w_1}\ldots\Omega_{w_m}\rangle =\sum_{\bar d}q^{\bar 
d}\langle\Omega_{w_1}\ldots\Omega_{w_m}\rangle_{\bar d}.
$$
Then $\wh{\s}_w(x ;q)$ is characterized by the condition
$$\langle\Omega_w\Omega_{w_1}\ldots\Omega_{w_m}\rangle 
=\langle\wh{\s}_w(x ,q)\Omega_{w_1}\ldots\Omega_{w_m}\rangle
$$
for any $w_1,\ldots ,w_m\in S_n$.

$\wh{\s}_w(x ;q)$ is called a geometric quantum Schubert polynomial. By 
definition $\wh{\s}_w(x;q)\in QH^{2l(w)}(\Fl )$.

\section{Proofs of Theorem 3 and quantum Cauchy formula.}

\begin{th} Let $I\in \t.$ Then 
\[ \langle \wt{e}_I(x),
\wt{\s}_w (x) \rangle _Q 
= \langle e_I(x), \s_w(x) \rangle , \] 
for any permutation $w\in S_n.$ 
\end{th} 
{\it Proof.} The proof is based on the arguments 
due to I.~Ciocan-Fontanine [C]. 
To begin with, let us remind his results. We consider 
the hyper-quot scheme 
${\cal HQ}_{\bar{d}}({\bf P}^1,Fl_n)$ 
associated to ${\bf P}^1$ with multidegree 
$\bar{d} = (d_1,\ldots,d_{n-1}).$ 
Let 
\[ {\bf C}^n \otimes {\cal O} \rightarrow 
{\cal T}_{n-1} \rightarrow \cdots \rightarrow 
{\cal T}_2 \rightarrow {\cal T}_1 \rightarrow 0 \] 
be the universal sequence of quotients on 
${\bf P}^1 \times {\cal HQ}_{\bar{d}}({\bf P}^1,Fl_n)$ 
and
$${\cal S}_i = {\rm Ker}\{ {\bf C}^n \otimes {\cal O} 
\rightarrow {\cal T}_{n-i} \} .
$$ 
We also consider the dual sequence 
\[ {\bf C}^n \otimes {\cal O} \rightarrow 
{\cal S}_{n-1}^{\ast} \rightarrow \cdots \rightarrow 
{\cal S}_1^{\ast}. \] 
We fix a flag 
\[ 0=V_0 \subset V_1 \subset \cdots \subset V_n= {\bf C}^n \] 
and define the subschema $D_w^{p,q}$ of 
${\bf P}^1 \times {\cal HQ}_{\bar{d}}({\bf P}^1,Fl_n)$ 
as the locus where ${\rm rank}(V_p \otimes {\cal O} \rightarrow 
{\cal S}_q^{\ast})\leq r_w(q,p),$ 
and $r_w(q,p):=\sharp \{ i \mid i\leq q, w_i\leq p \}.$ 
Let 
\[ D_w^{p,q}(t) = D_w^{p,q} \bigcap \{ \{t \} \times 
{\cal HQ}_{\bar{d}}({\bf P}^1,Fl_n) \} \] 
and 
\[ \overline{\Omega}_w(t) = \bigcap_{p,q=1}^{n-1}
D_w^{p,q}(t). \] 
Then the class of $\bar{\Omega}_w(t)$ in the Chow ring
$CH^{l(w)}({\cal HQ}_{\bar{d}}({\bf P}^1,Fl_n))$ is 
independent of $t\in {\bf P}^1$ and the flag 
$V_0\subset \cdots \subset V_n.$ 

The boundary ${\cal HQ}_{\bar{d}}({\bf P}^1,Fl_n)
\setminus M_{\overline{d}}({\bf P}^1,Fl_n)$ consists of 
$n-1$ divisors ${\bf D}_1,\ldots,{\bf D}_{n-1},$ 
which are birational to 
$${\bf P}^1 \times 
{\cal HQ}_{\bar{d}_1}({\bf P}^1,Fl_n),\ldots ,
{\bf P}^1 \times 
{\cal HQ}_{\bar{d}_{n-1}}({\bf P}^1,Fl_n)
$$ 
respectively, where $\bar{d}_{i}=(d_1,\ldots,d_i-1,\ldots,d_{n-1}).$ 
Let $x_i(t) = \overline{\Omega}_{s_i}(t) - 
\bar{\Omega}_{s_{i-1}}(t),$ 
then for any permutation $w\in S_n$ 
there exists an element 
$$G_w(t)\in CH_{\ast}(\bigcup_{i=1}^{n-2}{\bf D}_i)
$$ 
such that 
\[ \bar{\Omega}_{w^{-1}}(t)-\s_w(x_1(t),\ldots ,x_{n-1}(t)) 
= j_{\ast}(G_{w^{-1}}(t)), \]
where 
\[ j:\bigcup_{i=1}^{n-2}{\bf D}_i \rightarrow 
{\cal HQ}_{\bar{d}}({\bf P}^1,Fl_n) \] 
is the inclusion. 

Let $[m,k]\in S_n$ be the permutation 
\[ \left( 
\begin{array}{cccccccccccc}
1 & 2 & \ldots & m-k-1 & m-k & m-k+1 &\ldots & m & m+1 & \ldots 
 & n \\ 
1 & 2 & \ldots & m-k-1 & m & m-k &\ldots & m-1 & m+1 & \ldots 
 & n 
\end{array}
\right) . \] 
Then the geometric Schubert polynomial $\hat{\s}_{[m,k]^{-1}}(x)$ is 
the elementary symmetric function in $x_1,\ldots ,x_{m-1}$ 
of degree $k.$ Let $I=(i_1,\ldots,i_n).$ 
We have to calculate in $CH^{\ast}({\cal HQ}_{\bar{d}}
({\bf P}^1,Fl_n))$ 
\begin{equation} \left( \overline{  
\bigcap_{\nu =1}^{n-1}\Omega_{[n-\nu+1,i_{\nu}]^{-1}}} \right) (t)
\cdot \prod_{j=1}^{N}\overline{\Omega}_{w_j}(t_j) - 
\prod_{\nu =1}^{n-1}\s_{[n-\nu+1,i_{\nu}]}(x(t)) 
\cdot \prod_{j=1}^{N}\overline{\Omega}_{w_j}(t_j), \label{300}
\end{equation}
for distinct $t,t_1,\ldots, t_N$ and 
$w_1,\ldots ,w_N\in S_n$ such that 
$$\ds\sum_{\nu=1}^{n-1}i_{\nu}
+\ds\sum_{j=1}^{N}l(w_j) = n(n-1)/2 + 2\ds\sum_{k=1}^{n-1} d_k,
$$ 
where $\cap$ is the classical intersection product and 
\[ \left( \overline{  
\bigcap_{\nu =1}^{n-1}\Omega_{[n-\nu+1,i_{\nu}]^{-1}} } \right) (t) \] 
is the corresponding degeneracy locus on 
${\cal HQ}_{\bar{d}}({\bf P}^1,Fl_n).$ 
In order to calculate the expression (\ref{300}), we are going to prove
at first that
\[ \left( \overline{ \bigcap_{\nu =1}^{n-1}
\Omega_{[n-\nu +1,i_{\nu}]^{-1}}} \right) (t) \cdot
\prod_{j=1}^{N}\bar{\Omega}_{w_j}(t_j) -
\prod_{\nu =1}^{n-1}
\bar{\Omega}_{[n-\nu +1,i_{\nu}]^{-1}} (t)
\cdot
\prod_{j=1}^{N}\bar{\Omega}_{w_j} (t_j) = 0 .\]
The LHS of the last expression can be computed as the number of points in
\[ \prod_{\nu =1}^{n-1}
\bar{\Omega}_{[n-\nu +1,i_{\nu}]^{-1}}(t)
\cdot
\prod_{j=1}^{N}\bar{\Omega}_{w_j} (t_j) \]
suppoted on $\cup_{i=1}^{n-2}{\bf D}_i.$ Let $j_{\bar{d}_k}$ be
the natural rational map
\[ {\bf P}^1 \times {\cal HQ}_{\bar{d}_k}
({\bf P}^1,Fl_n) - \rightarrow
{\cal HQ}_{\bar{d}}({\bf P}^1,Fl_n). \]
From Remark 3 in [C],
\[ j_{\bar{d}_k}^{-1}(\bar{\Omega}_{[n-\nu+1,i_{\nu}]^{-1}}(t)) = \]
\[\left\{
\begin{array}{cc}
{\bf P}^1 \times \bar{\Omega}_{[n-\nu+1,i_{\nu}]^{-1}}(t) \cup
\bigcap_{n-\nu-i_\nu+1\leq p \leq n-\nu-1}
D_{[n-\nu+1,i_{\nu}]^{-1}}^{p,n-\nu-1}(t), & {\rm if} \;
k=n-\nu,
\\
{\bf P}^1 \times \bar{\Omega}_{[n-\nu+1,i_{\nu}]^{-1}}(t) &
{\rm otherwise,}
\end{array}
\right.  \]
where 
\[ \bigcap_{n-\nu-i_\nu+1\leq p \leq n-\nu-1}
D_{[n-\nu+1,i_{\nu}]^{-1}}^{p,n-\nu-1}(t) = 
\{ t\} \times \bar{\Omega}_{[n-\nu,i_{\nu}-1]^{-1}}(t) .\]
Because, by assumption, 
\[\sum_{\nu =1}^{n-1}i_{\nu}+\sum_{\nu =1}^{n-1}l(w_j) 
= \frac{n(n-1)}{2} + 2 \sum_{k=1}^{n-1}d_k, \]
we have 
\[ \bar{\Omega}_{[k,i_{n-k}-1]^{-1}}(t)\cdot 
\prod_{\nu \not= n-k}
\bar{\Omega}_{[n-\nu +1,i_{\nu}]^{-1}}(t)
\cdot
\prod_{j=1}^{N}\bar{\Omega}_{w_j} (t_j)
=0 \] 
on  the hyper-quot scheme ${\cal HQ}_{\bar{d}_k}({\bf P}^1,Fl_n).$ 
Hence, we have equality
\[ \left( \overline{ \bigcap_{\nu =1}^{n-1}
\Omega_{[n-\nu +1,i_{\nu}]^{-1}}} \right) (t) \cdot
\prod_{j=1}^{N}\bar{\Omega}_{w_j}(t_j) =
\prod_{\nu =1}^{n-1}
\bar{\Omega}_{[n-\nu +1,i_{\nu}]^{-1}} (t)
\cdot
\prod_{j=1}^{N}\bar{\Omega}_{w_j} (t_j) . \]

Our next observation is that the intersection number in the RHS of 
the last equality is equal to
$$\prod_{\nu =1}^{n-1}\overline\Omega_{[n-\nu +1,i_{\nu}]^{-1}}(s_{\nu})
\cdot\prod_{j=1}^N\overline\Omega_{w_j}(t_j),
$$
where we can chose $s_1,\ldots ,s_{n-1}$, $t_1,\ldots ,t_N\in P^1$ to be
the pairwise distinct points, since the class $[\Omega_w(t')]$ in the
Chow ring $CH^*({\cal H}Q_{\overline d}(P^1,Fl_n))$ does not depends on 
the chose of $t'\in P^1$. 

Now we are going to use the following identity
$$\prod_{k=1}^ma_k-\prod_{k=1}^mb_k=\sum_{k=1}^m\prod_{j=1}^{k-1}b_j
(a_k-b_k)\prod_{j=k+1}^ma_j.
$$
Let us take in the last equality $m=n-1$
\begin{eqnarray*}
&&a_k:=\overline\Omega_{[n-k+1, i_k]^{-1}}(s_k),\\
&&b_k:=\s_{[n-k+1,i_k]}(x(s_k)).
\end{eqnarray*}
Then we obtain the following equality
\begin{eqnarray*}
& & \prod_{\nu =1}^{n-1}\overline\Omega_{[n-\nu +1,i_{\nu}]^{-1}}(s_{\nu})
\prod_{j=1}^N\overline\Omega_{w_j}(t_j)-\prod_{\nu =1}^{n-1}
\s_{[n-\nu +1,i_{\nu}]}(x(s_{\nu}))\prod_{j=1}^N\overline\Omega_{w_j}(t_j)\\
&=&\left\{\sum_{k=1}^{n-1}\prod_{j=1}^{k-1}\overline\Omega_{[n-j+1,i_j]^{-1}}
(s_j)\cdot j_*(G_{[n,i_k]^{-1}}(t))\prod_{j=k+1}^{n-1}\s_{[n-j+1,i_j]}
(x(s_j))\right\}\cdot B,
\end{eqnarray*}
where $B:=\prod_{j=1}^N\overline\Omega_{w_j}(t_j)$.

The contributions from $j_*((G_{[n,i_k]^{-1}}(t))$, $1\le k\le n-1$, can
be computed by using the arguments in [C]. Inded, as in the proof of 
Theorem~4 in [C], the intersection number ($1\le k\le n-1$)
$$\prod_{l=1}^{k-1}\overline\Omega_{[n-l+1,i_l]^{-1}}
(s_l)\cdot j_*(G_{[n,i_k]^{-1}}(t))\prod_{l=k+1}^{n-1}\s_{[n-l+1,i_l]}
(x(s_l))\cdot B
$$
is the number of points in
$$\prod_{l=1}^{k-1}\overline\Omega_{[n-l+1,i_l]^{-1}}
(s_j)\prod_{l=k}^{n-1}\s_{[n-l+1,i_l]}(x(s_l))\cdot B
$$
supported on $\bigcup_{j=1}^{n-2}D_i$. Hence, by induction, we have the
following identity for correlation functions
$$\langle\left(\bigcap_{\nu =1}^{n-1}\Omega_{[n-\nu +1,i_{\nu}]^{-1}}\right)
\cdot\prod_{j=1}^N\overline\Omega_{w_j}\rangle =\langle\prod_{\nu =1}^{n-1}
\wt{\s}_{[n-\nu +1,i_{\nu}]}\cdot\prod_{j=1}^N\overline\Omega_{w_j}\rangle ,
$$
where $\bigcap$ is the classical intersection product and $\cdot$ is the 
product in 
$${\rm Sym}~H^*(Fl_n,{\bf Z})[q_1,\ldots ,q_{n-1}]. 
$$
The last equality
for correlation functions is equivalent to the following one
\[ m_q(\Omega_{[n,i_1]^{-1}}\cap \Omega_{[n-1,i_2]^{-1}}, \cap\ldots\cap
\Omega_{[l,i_{n-1}]^{-1}}, *) =m_q(\wt e_{i_1}\cdot\wt e_{i_2}\cdots
\wt e_{i_{n-1}},*).
\] 
This completes the proof.
\qed \medskip \\ 
{\bf Remark 10} \ Using the similar geometrical arguments we can prove an 
analog of Theorem~9 for the quantum complete polynomials: \smallskip \\
{\bf Theorem 9'} {\it Let $I\subset\delta_n$. Then
$$\langle\wt h_I(x),\wt{\s}_w(x)\rangle_Q=\langle h_I(x),\s_w(x)\rangle,
$$
for any permutation $\omega\in S_n$.}

It is easy to see that Theorem 5' is a corollary of Theorem 9'.

\section{Correlation functions.}

\subsection{Higher genus correlation function and the 
Vafa--Intriligator type formula.}

Fix a Riemann surface $C$ of genus $g$. We denote by $M_d(C,F)$ the moduli 
space of morphism from $C$ to $\Fl $. One can define the higher genus 
Gromov--Witten invariants by method which is similar to that in the 
case of genus zero, [RT].

We have the following recursion relation for higher genus correlation 
function corresponding to the generating function for higher genus 
Gromov--Witten invariants
$$\langle\Omega_{w_1}\ldots\Omega_{w_N}\rangle_g=
\sum_{v\in 
S_n}\langle\Omega_{w_1}\ldots\Omega_{w_N}\Omega_v\Omega^*_v\rangle_{g-1}
$$
(cf. Ruan--Tian [RT]).

From Corollary \ref{107} and Theorem~\ref{011} 
we can deduce the Vafa--Intriligator type formula 
for higher genus correlation functions, namely,
let $\langle P(x_1,\ldots ,x_n)\rangle_g$ be the genus $g$
correlation function corresponding to a polynomial $P$, then

\vskip 0.2cm
\hskip 2.0cm $\langle P(x_1,\ldots ,x_n)\rangle_g=\Res_{\wt I}(P\Phi^g )=$
$$=\sum_{\wt e_1=\cdots =\wt e_n=0}P(x_1,\ldots ,x_n)\det\left 
(\frac{\partial\wt e_i}{\partial x_j}\right)^{-1}\left(\Phi (x_1,\ldots 
,x_n)\right)^g,
$$
where $\Phi 
(x)=\langle \wt{\s}_{w_0}(x,y),\wt{\s}_{w_0}(x,y)\rangle^{(y)}=
\ds\sum_{w\in S_n}\wt{\s}_w(x)\wt{\s}_{w_0w}(x)=C^{(q,q)}(x,x)$.

To simplify the formula above, we use the following observations

$\bullet$ $\det\left(\ds\frac{\partial\wt e_i}{\partial x_j}\right)\equiv
n!\:\wt{\s}_{w_0}(x)~(\mod\wt I);$

$\bullet$ $\Phi (x):=C^{(q,q)}(x,x)\equiv n!\:\wt{\s}_{w_0}(x)~(\mod\wt I).$

Hence, we obtain
\begin{th} (Higher genus Vafa-Intriligator formula)
$$\langle P(x_1,\ldots ,x_n)\rangle_g=\Res_{\wt I}(P\Phi^g)=
\sum_{\wt e_1=\cdots =\wt e_n=0}P(x_1,\ldots ,x_n)\left(\wt{\s}_{w_0}(x)
\right)^{g-1},
$$
where $\wt{\s}_{w_0}(x)=\wt e_{\delta}(x):=\wt e_1(x_1)\wt e_2(x_1,x_2)
\ldots\wt e_{n-1}(x_1,\ldots ,x_{n-1})$, and $\wt e_i(z)$ is the quantum
elementary polynomial of degree $i$ in the variable $z=(z_1,\ldots ,z_m)$,
see Section 5.2.
\end{th}
{\bf Remark 11} \ The polynomial
$$C^{(q,q')}(x,y):=\sum_{w\in S_n}\wt{\s}_w^{(q)}(x)
\wt{\s}_{w_0w}^{(q')}(y)
$$
corresponds to the dual class of the diagonal in the quantum cohomology 
ring $QH^*(\Fl \times\Fl ,(q,q'))=QH^*(\Fl ,q)\otimes QH^*(\Fl ,q')$.

\subsection{Witten--Dijkgraaf--Verlinde--Verlinde equations \break 
for symmetric group.}

The Witten--Dijkgraaf--Verlinde--Verlinde equations (WDVV-equations) are 
equations on the correlation functions 
$\langle\wt{\s}_u\wt{\s}_v\wt{\s}_w\rangle\in{\bf Z}[q_1,\ldots 
,q_{n-1}]$, where $u,v,w\in S_n$. The correlation functions satisfy the 
following conditions

1) Normalization: 
$$\langle 1\wt{\s}_v\wt{\s}_w\rangle =\langle 
\s_v,\s_w\rangle .
$$

2) Initial data: 
$$\langle \wt{\s}_{s_k}\wt{\s}_{s_k}\wt{\s}_{w_0}\rangle =q_k.
$$

3) Degree conditions:
$$\langle\wt{\s}_u\wt{\s}_v\wt{\s}_w\rangle =0
$$
if either $l(u)+l(v)+l(w)<l(w_0)$, or difference $l(u)+l(v)+l(w)-l(w_0)$
is an odd positive integer.

4) WDVV-equations: 
$$\sum_v\langle\wt{\s}_{w_1}\wt{\s}_{w_2}\wt{\s}_v\rangle
\langle\wt{\s}_{w_0v}\wt{\s}_{w_3}\wt{\s}_{w_4}\rangle =
\sum_v\langle\wt{\s}_{w_2}\wt{\s}_{w_3}\wt{\s}_v\rangle
\langle\wt{\s}_{w_0v}\wt{\s}_{w_1}\wt{\s}_{w_4}\rangle .
$$
for any $w_1,w_2,w_3,w_4\in S_n$. \smallskip \\ 
{\bf Conjecture 1} \ {\it Conditions 1)--4) uniquely determine the correlation 
functions $\langle\wt{\s}_u\wt{\s}_v\wt{\s}_w\rangle$.}
\smallskip \\
{\bf Remark 12} \ 1) Correlation function 
$\langle\wt{\s}_{w_1}\wt{\s}_{w_2}\wt{\s}_{w_3}\rangle $ is a 
generating function for the Gromov--Witten invariants (see Definition 4):
$$\langle\wt{\s}_{w_1}\wt{\s}_{w_2}\wt{\s}_{w_3}\rangle :=
\ds\sum_{\overline d}q^{\overline d}\langle\wt{\s}_{w_1}\wt{\s}_w\wt{\s}_{w_3}
\rangle_{\overline d}.
$$
2) More generally,
$$\langle\wt{\s}_{s_k}\wt{\s}_{s_{ij}}\wt{\s}_{w_0}\rangle 
=\left\{\begin{array}{ll} q_i\ldots q_{j-1}, & {\rm if} \ 1\le i\le 
k<j\le n,\\ 0, & {\rm otherwise}.\end{array}\right.
$$

\subsection{Residue formula.} 

\begin{th}\label{011} Correlation function 
$\langle P(x_1,\ldots ,x_n) \rangle$ is given by 
the formula 
\[ \langle P(x_1,\ldots ,x_n) \rangle = 
\sum_{{\wt e}_1(p)=\cdots ={\wt e}_{n}(p)=0}{\rm Res}_p 
\left( \frac{ P(x_1,\ldots ,x_n)
dx_1 \wedge \cdots \wedge dx_n }
{{\wt e}_1 \cdots {\wt e}_n} \right). \]
\end{th} 
{\it Proof.} If the polynomial $P(x_1,\ldots ,x_n)$ 
is in the ideal generated by ${\wt e}_1,\ldots ,{\wt e}_n$, 
then the left and right hand sides of the formula are
zero. Hence, it is enough to prove that
$$\sum{\rm Res}_p\left(\frac{x_1^{\nu_1}\cdots x_{n-1}^{\nu_{n-1}}dx_1
\wedge\cdots\wedge dx_n}{\wt e_1\cdots \wt e_n}\right)
$$
\[ = \left\{
\begin{array}{ll}
1, & {\rm if}\; \;(\nu_1,\ldots ,\nu_{n-1})
= (n-1,n-2,\ldots,1) \\ 
0, & {\rm if}\; \; 0\le \nu_1+\cdots + \nu_{n-1} < 
n(n-1)/2, \; \; 0\leq \nu_i \leq n-i.
\end{array} 
\right. \] 
We can extend the meromorphic form 
\[ \omega = \frac{ x_1^{\nu_1}\cdots x_{n-1}^{\nu_{n-1}}
dx_1 \wedge \cdots \wedge dx_n }
{{\wt e}_1 \cdots {\wt e}_n} \] 
on the affine space ${\bf A}^{n}_{(x_1,\ldots ,x_n)}$ 
to $({\bf P}^1)^n.$ 
For each subset $I \subset \{ 1,\ldots,n \} ,$ 
we consider the coordinate chart $U_J = 
{\bf A}^{n}_{(z_1^J,\ldots ,z_n^J)},$ where 
\[ z_i^J = \left\{ 
\begin{array}{ll} 
x_i, & {\rm if}\; \; i \not\in J, \\ 
1/x_i, & {\rm otherwise.} 
\end{array} 
\right. \] 
Then 
\[ ({\bf P}^1)^n = \bigcup_{J} U_J. \] 
Let 
\[ \bar{e}_j^J(z_1^J,\ldots ,z_n^J) = 
(\prod_{i\in J} z_i^J) {\wt e}_j(x_1,\ldots ,x_n) . \] 
Then 
\[ \omega = (-1)^{\sharp J} 
\frac{x_1^{\nu_1}\cdots x_{n-1}^{\nu_{n-1}}
\ds(\prod_{i\in J}z_i^J)^{n-2}dz_1^J
\wedge \cdots \wedge dz_n^J}
{\bar{e}_1^J \cdots \bar{e}_n^J } \] 
on $U_{\phi}\cap U_J.$ 
If $\sharp J = j,$ then there exists a polynomial 
$Q_i(z_1^J,\ldots ,z_n^J)$ for $i\in J$ 
such that 
\[ \bar{e}_j^I(z_1^J,\ldots ,z_n^J) = 
1+\sum_{i\in J}z_i^J Q_i. \] 
This follows from 
\[ {\wt e}_j(x_1,\ldots ,x_n) = e_j(x_1,\ldots ,x_n) +
({\rm terms} \; \; {\rm of} \;\; {\rm lower} \; \; 
{\rm degree}). \] 
Therefore $\bar{e}_1^J,\ldots \bar{e}_n^J$ 
do not have common zero on 
\[ B_J = \{ (z_1^J,\ldots ,z_n^J) \in U_J \mid 
 z_i^J = 0 , \; i\in J \}. \] 
From the residue theorem, if 
$0\leq \nu_1+\cdots + \nu_n < n(n-1)/2,$ 
$0\leq \nu_i \leq n-i,$ then 
\[ \sum {\rm Res}_p 
\left( \frac{ x_1^{\nu_1}\cdots x_{n-1}^{\nu_{n-1}}
dx_1 \wedge \cdots \wedge dx_n }
{{\wt e}_1 \cdots {\wt e}_n} \right)=0. \] 
On the other hand, 
\[ \sum_{{\wt e}_i(p)=0} {\rm Res}_p 
\left( \frac{ x_1^{n-1}\cdots x_{n-1}
dx_1 \wedge \cdots \wedge dx_n}{{\wt e}_1\cdots {\wt e}_n} \right) = 
-\sum_p{\rm Res}_p \omega , 
\] 
where $p$ runs over the common zeros of 
$\bar{e}^J_1, \ldots ,\bar{e}^J_n$ in 
$\ds\bigcup_{1\in J}B_J.$ 
Let $y_1 = 1/x_1,$ $z=(y_1,x_2,\ldots ,x_n)$ and 
\[ \bar{e}^{*}_1(z) = 
y_1(1+ y_1(x_2+ \cdots + x_n)). \] 
Then we have \ \ 
$ -\sum_p{\rm Res}_p \omega = $
\[ \sum_{\scriptsize{\begin{array}{cc}\bar{e}^{\ast}_1 = 
\bar{e}_2^{\{ 1\} } 
= \cdots =\bar{e}_n^{\{ 1\} } = 0 \\ {\rm in} 
\; \; {\rm the }\; \; {\rm locus} \; \; 
\{ y_1 =0 \}\end{array}} } 
{\rm Res}_p \left(\frac{ x_2^{n-2}\cdots x_{n-1}
dy_1 \wedge \cdots \wedge dx_n }
{\bar{e}^{\ast}_1(z)\cdot
 \bar{e}^{\{ 1\} }_2(z)\cdots 
 \bar{e}^{\{ 1\} }_N(z)}
\right) \] 
\[ = \sum_p{\rm Res}_p
\left( \frac{ x_2^{n-2}\cdots x_{n-1}
dx_2 \wedge \cdots \wedge dx_n }
{{\wt e}_1(x_2,\ldots ,x_n) \cdots 
{\wt e}_{n-1}(x_2,\ldots ,x_n)} 
\right) = 1, \] 
by induction. \qed \bigskip 

From the residue formula, the correlation function is 
given by the quantum residue ${\rm Res}_{\wt I},$ namely, 
\[ \langle P(x_1, \ldots , x_{n-1}) \rangle = {\rm Res}_{\wt I} 
P(t_1, \ldots ,t_{n-1}) . \] 
In order to relate the quantum residue with the 
classical one, we consider the quantum residue generating 
function 
\[ \Psi(t) = \langle \prod_{i=1}^{n-1} \frac{t_i}
{t_i - x_i} \rangle 
= \sum_{\nu \in ({\bf Z}_{\geq 0})^{n-1}}
\langle x^{\nu} \rangle t^{-\nu}.\] 
Then, we have 
\[ {\rm Res}_{\wt I} P(t_1,\ldots ,t_{n-1}) = 
{\rm Res}_I \left( P(x_1,\ldots ,x_{n-1}) 
\Psi(x) \right) . \] 
Hence, it is important to determine the generating function 
$\Psi(t).$ Let 
\[ f_i(t) = t^n + 
\sum_{j=0}^{n-1} \gamma_{n-j}^{(i)} t^j \] 
be the characteristic polynomial of the 
quantum multiplication by $x_i$ with respect to the basis 
consisting of the quantum Schubert polynomials. 
Let us consider the $(n!+1)\times (n!+1)$-matrix $C_n(t)$ 
such that 
\begin{eqnarray*}(C_n(t))_{1,j} &=& \frac{(-1)^{j-1} t^{n-j+2}}{(j-1)!}; \\
\\ 
(C_n(t))_{i,j} &=& \left\{  
\begin{array}{ll}
\ds\frac{(-1)^n}{(j-1)!} 
{i-2 \choose n-j+1} t^{j-2}, & {\rm if} 
\; \; i\geq 2, \; i+j \geq n+2,  
 \\ \\
0, & {\rm otherwise.} 
\end{array} 
\right. 
\end{eqnarray*} 
We define the differential operator $D_i$ by 
\[ D_i = (\gamma_{n!}^{(i)},\gamma_{n!-1}^{(i)},\ldots,1)
\cdot C_n(t_i) \cdot 
 ^{t}\! (1,\partial / \partial t_i,\ldots 
,(\partial / \partial t_i )^{n!} ). \] 
\begin{pr} The generating function 
$\Psi(t)$  
satisfies the system 
of differential equations 
\[ D_i \Psi(t) = 0, 
\; \; 1\leq i \leq n-1. \] 
Conversely, these differential equations and the initial 
values $\langle x_1^{\nu_1}\cdots x_{n-1}^{\nu_{n-1}} \rangle $ 
for $0\leq \nu_i \leq n!-1$ 
determine the generating function uniquely.  
\end{pr} 
{\it Proof.} Let $x$ be a variable. Since 
\[ \left( \begin{array}{c}
1 \\ x \\ x^2 \\ x^3 \\ \vdots \\ x^{n!} 
\end{array} \right) 
= C_n \cdot 
\left( \begin{array}{c}
t(t-x)^{n!} \\ -x(t-x)^{n!-1} \\ 2! x(t-x)^{n!-2} 
\\ -3! x(t-x)^{n!-3} \\ \vdots \\ (-1)^{n!}(n!)!x
\end{array} \right) , \] 
we have 
\[ D_i \Psi(t)
= \langle \frac{f_i(x_i)}{(t_i - x_i)^{n!+1}}
\prod_{j \not= i}\frac{t_j}{t_j - x_j} 
 \rangle 
= 0. \] 
On the other hand, the recursive relations 
\[ \langle f_i(x_i) P(x_1,\ldots ,x_{n-1}) 
\rangle = 0 \] 
with the initial values 
$\langle x_1^{\nu_1}\cdots x_{n-1}^{\nu_{n-1}} \rangle $ 
for $0\leq \nu_i \leq n!-1$ 
determine the correlation function uniquely. \qed
\smallskip \\ 
{\bf Remark 13} \ We can also consider another generating 
function 
\[ \langle \exp(x_1 t_1 + \cdots + x_n t_n) \rangle .\] 
This is the generating volume function in [GK]. This 
generating function satisfies 
\[ {\wt e}_i \left( \frac{\partial}{\partial t_1},\cdots , 
\frac{\partial}{\partial t_n}\right) 
\langle \exp (x_1 t_1 + \cdots + x_n t_n) \rangle 
= 0, \] 
for $1\leq i \leq n.$ \bigskip 
 
In the case of $n=3,$ we can calculate the generating 
function $\Psi (t)$ explicitly. The results of calculation one can find
in the Appendix B.

\section{Extended Ehresman--Bruhat order and\break quantum Pieri rule.}

Let us remind that the Ehresman--Bruhat order denoted by $\le$, 
is the partial order on $S_n$ 
that is the transitive closure of the relation $\rightarrow$. Relation 
$v\to w$ means that

1) $l(w)=l(v)+1$,

2) $w=v\cdot t$ where $t$ is a transposition.

In other words, if $v$ and $w$ are permutations, $v\le w$ means that 
there exists $r\ge 0$ and $v_0,v_1,\ldots ,v_r$ in $S_n$ such that
$$v=v_0\to v_1\to\cdots\to v_r=w.
$$

Now let us define the extended Ehresman--Bruhat order $v\Leftarrow w$ on 
$S_n$. First of all, we define a relation $v\leftarrow w$ (see, also, [FGP]). 
Relation $v\leftarrow w$ means that

1) $w=v\cdot t$, where $t$ is a transposition,

2) $l(w)\ge l(v)+l(v^{-1}w)$.
\smallskip \\
{\bf Remark 14} \
$i)$ It follows from [M], (1.10), that condition 2) is 
equivalent to the following one \smallskip \\
$2'$) $w(i)<w(j)$ and for all $k$ such that $i<k<j$ 
we have $w(i)<w(k)<w(j)$. 

$ii)$ If $w=vt_{i,i+1}$ and $l(w)=l(v)+1$ (i.e. $v\to w$ in the Bruhat 
order), then we have also an arrow $v\leftarrow w$. This is clear because 
in our case we have $l(w)=l(v)+l(t_{i,i+1})$.

\medskip 
{\bf Example.} Symmetric group $S_3$ (see Figure 1).

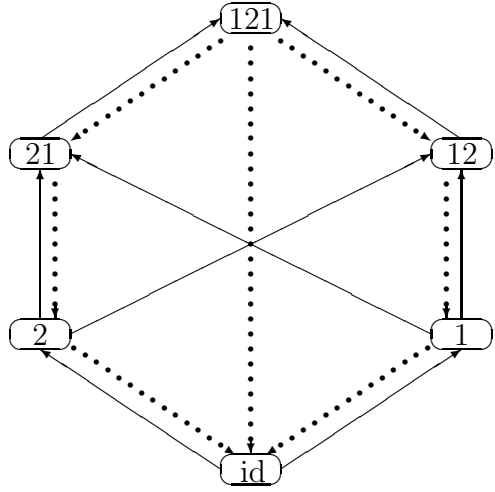
\begin{figure}[hbtp]
\setlength{\unitlength}{0.4cm}
\begin{picture}(16,16)(-7.5,1)
\put(7,1){\vector(-3,2){6}}
\put(1,6){\vector(0,1){5}}
\put(1,12){\vector(3,2){6}}
\put(15,6){\vector(0,1){5}}
\put(15,12){\vector(-3,2){6}}
\put(9,1){\vector(3,2){6}}
\put(2,5.5){\vector(2,1){12}}
\put(14,5.5){\vector(-2,1){12}}
\put(8,1){\oval(2,1)}
\put(7.65,0.6){\hbox{\rm id}}
\multiput(8,2)(0,0.5){27}{\circle*{0.2}}
\put(8,2){\vector(0,-1){0.5}}
\multiput(7,1.8)(-0.375,0.25){14}{\circle*{0.2}}
\put(7,1.8){\vector(3,-2){0.5}}
\multiput(9,1.8)(0.375,0.25){14}{\circle*{0.2}}
\put(9,1.8){\vector(-3,-2){0.5}}
\put(1,5.5){\oval(2,1)}
\put(0.75,5.15){\hbox{2}}
\multiput(1.5,6.5)(0,0.5){9}{\circle*{0.2}}
\put(1.5,6.5){\vector(0,-1){0.5}}
\put(1,11.5){\oval(2,1)}
\put(0.5,11.15){\hbox{21}}
\put(8,16){\oval(2,1)}
\put(7.25,15.65){\hbox{121}}
\multiput(2.5,12.25)(0.375,0.25){13}{\circle*{0.2}}
\put(2.5,12.25){\vector(-3,-2){0.5}}
\multiput(13.5,12.25)(-0.375,0.25){13}{\circle*{0.2}}
\put(13.5,12.25){\vector(3,-2){0.5}}
\put(15,11.5){\oval(2,1)}
\put(14.5,11.15){\hbox{12}}
\put(15,5.5){\oval(2,1)}
\put(14.75,5.15){\hbox{1}}
\multiput(14.5,6.5)(0,0.5){9}{\circle*{0.2}}
\put(14.5,6.5){\vector(0,-1){0.5}}

\end{picture}
\medskip
\caption{Extended Ehresman-Bruhat order for $S_3$. }
\end{figure}

We define a weight of an arrow $v\leftarrow w$, denoted by 
$wt(v\leftarrow w)$, to be equal to the product 
$q_i\ldots q_{i+s-1}$, if $t=t_{ij}$ and $2s:=l(w)+1-l(v)$. We assume 
that weight of any arrow $v\to w$ is equal to 1 (see, also, [FGP]).

Let us say that an arrow $v\leftarrow w$ (resp. $v\to w$) 
has  a color $k$ if $w=vt_{ij}$ and $1\le i\le k<j\le n$. 

Extended Ehresman--Bruhat order on $S_n$ (notation 
$v\Leftarrow w$)
is the transitive closure of the 
relations $\leftarrow$, and $\rightarrow$. In other words, there 
exists $r\ge 0$ and 
$v_0,v_1,\ldots ,v_r$ in $S_n$ such that

\begin{equation} v=v_0\rightleftharpoons v_1\rightleftharpoons 
v_2\rightleftharpoons\cdots\rightleftharpoons v_r=w,
\label{1}
\end{equation}
where symbol $v_i\rightleftharpoons v_{i+1}$ means either $v_i\rightarrow 
v_{i+1}$ or $v_i\leftarrow v_{i+1}$. 

For given pair $v\Leftarrow w$, we consider a sequence of arrows (\ref{1})
as a path between $v$ and $w$ (notation $v\mapsto w$) in the extended
Ehresman--Bruhat order and call it as a BE--path (Bruhat--Ehresman path).
We denote the number $r$ in a 
representation (\ref{1}) by $l(v\mapsto w$).

Let us define a weight of a BE--path $v\mapsto w$ as follows 
$$wt(v\mapsto w)=
\ds\prod_{i=0}^{r-1}wt(v_i\rightleftharpoons v_{i+1}).
$$ 
We will say that 
BE--path $v\mapsto w$ has a color $k$, notation $v~{\buildrel k
\over\mapsto }~w$, if in the representation (\ref{1}) 
all arrows $v_i\rightleftharpoons v_{i+1}$ $(i=0,\ldots ,r-1)$ have the 
same color $k$.

\begin{th} (Quantum Pieri's rule). Let us consider the Grassmanian permutation 
$[b,d]=(1,2,\ldots ,b-d-1,b,b-d,b-d+1\ldots ,b-1,b+1,\ldots ,n)$, 
for $2\leq b\leq n,$ $1\leq d\leq b.$ Then
$$\wt{\s}_{[b,d]}\cdot\wt{\s}_v\equiv\sum wt(v~{\buildrel b
\over\mapsto }~w )\wt{\s}_w~(\mod \wt I_n),
$$
where the sum runs over all BE--paths $v~{\buildrel b\over\mapsto }~w$, 
s.t.

1) $l(v\mapsto w)=d$;
\vskip 0.2cm

2) if $v_l=v_{l+1}(i_lj_l)$ ($l=0,\ldots ,d-1$), then all $i_l$ are 
different.
\end{th}
(Note that $\s_{[b,d]}=e_d(x_1,\ldots ,x_{b-1})$).

{\it Sketch of the proof.} It is enough to consider the case $d=1$ 
(induction!). In the 
case $d=1$, we use a quantum analog of Kohnert--Veigneau's method [KV]. 
Namely, at first we prove the quantum Pieri rule (for $d=1$) for 
double quantum Schubert polynomials and then take $y=0$ (see Theorem~4).

\begin{pr} (Quantum Pieri's rule for $\wt{\s}_{w_0}(x,y)$).
\begin{equation} (x_j+y_{n+1-j})\wt{\s}_{w_0}(x,y)\equiv
\sum_{i<j}q_{ij}\wt{\s}_{w_0t_{ij}}(x,y)-\sum_{j<k}q_{jk}
\wt{\s}_{w_0t_{jk}}(x,y)\ (\mod \wt J), \label{200}
\end{equation}
where $q_{ij}:=q_iq_{i+1}\ldots q_{j-1}$, if $i<j$;\\ 
$\wt J$ is the ideal in the ring 
${\bf Z}[x_1,x_2,\ldots ,x_n,y_1,\ldots ,y_n,q_1,\ldots ,q_{n-1}]$ 
generated by
$$e_i(x_1,\ldots x_n~|~q_1,\ldots ,q_{n-1})+(-1)^{i-1}e_i(y_1,\ldots ,y_n),
\ 1\le i\le n,
$$
and $e_k(x_1,\ldots ,x_n~|~q_1,\ldots ,q_{n-1})$ is the $k$-th quantum 
elementary symmetric function.
\end{pr}

Applying to (\ref{200}) the generalized Monk formula (see Section 2.2), 
we obtain
\begin{col} (Equivariant quantum Pieri's rule)
\begin{eqnarray*}& &x_j\wt{\s}_w(x,y)+y_{w_j}\wt{\s}_w(x,y)\equiv\\
& &\equiv \sum_{j<k,~l(wt_{jk})=l(w)+1}\wt{\s}_{wt_{jk}}(x,y)+
\sum_{j<k,~l(w)=l(wt_{jk})+l(t_{jk})}q_{jk}\wt{\s}_{wt_{jk}}(x,y) \\
& &-\sum_{i<j,~l(wt_{ij})=l(w)+1}\wt{\s}_{wt_{ij}}(x,y)-
\sum_{i<j,~l(w)=l(wt_{ij})+l(t_{ij})}q_{ij}\wt{\s}_{wt_{ij}}(x,y)\ 
(\mod \wt J).
\end{eqnarray*}
\end{col}
\vskip 0.3cm
{\bf Remark 15} \  
$i)$ $\wt{\s}_{[b,d]}=e_d(x_1,\ldots ,x_b~|~q_1,\ldots ,q_{b-1})$ 
coincides with quantum elementary symmetric function.

$ii)$ It is clear that (initial data) 
$$\wt{\s}_{sk}^2=\wt{\s}_{s_{k+1}s_k}+\wt{\s}_{s_{k-1}s_k}+q_k,
$$ 
i.e. $\langle\wt{\s}_k\wt{\s}_k\wt{\s}_{w_0}\rangle =q_k$.

$iii)$ To our knowledge, in the classical case $q=0$, the Pieri rule for
Schubert polynomials was first stated in [LS1], (2.2). Our formulation
of Theorem~12 is very close to that given in [BB]. The difference is:
we use the paths in the extended Ehresman--Bruhat order (quantum case)
instead of the paths in the ordinary Ehresman--Bruhat order (classical
case). Very transparent proof of Monk's formula one can find in the
I.~Macdonald book [M1], (4.15). It is the proof that was generalized
in [FGP] to the case of quantum Schubert polynomials. Recently, F.~Sotile
[S] gave a proof of the Pieri rule based on geometrical approach.

\vskip 0.3cm
{\it Acknowledgement.} We would like to acknowledge our special 
indebtedness to Dr.~N.A.~Liskova for the inestimable help in 
preparing the manuscript for publication. 

\newpage
\section*{Appendix A} 
{\large {\bf Quantum double Schubert polynomials for} $S_4.$ }
\begin{eqnarray*}
\wt{\s}_{121321}(x,y) & = & \Delta_1(y_3\mid x_1) \Delta_2(y_2\mid x_1,x_2) 
\Delta_3(y_1\mid x_1,x_2,x_3), \\ 
\wt{\s}_{21321}(x,y) & = & q_1^2 x_1 + q_1 q_2 x_1 - 
 q_2 x_1^3 + 2 q_1 x_1^2 x_2 + 
x_1^3 x_2^2 + q_1 x_1^2 y_1 - q_2 x_1^2 y_1  \\ 
 &+ &   q_1 x_1 x_2 y_1 + x_1^3 x_2 y_1 + 
 x_1^2 x_2^2 y_1 + q_1 x_1 y_1^2 + x_1^2 x_2 y_1^2 + 
 q_1 x_1^2 y_2 \\ 
 &- &  q_2 x_1^2 y_2 + 
 q_1 x_1 x_2 y_2 +  x_1^3 x_2 y_2 + x_1^2 x_2^2 y_2 - 
  q_2 x_1 y_1 y_2 + x_1^3 y_1 y_2\\ 
 &+ & 2 x_1^2 x_2 y_1 y_2 +  x_1 x_2^2 y_1 y_2 + x_1^2 y_1^2 y_2 + 
  x_1 x_2 y_1^2 y_2 + q_1 x_1 y_2^2\\ 
 &+ &  x_1^2 x_2 y_2^2 + x_1^2 y_1 y_2^2 + x_1 x_2 y_1 y_2^2 + 
 x_1 y_1^2 y_2^2 + q_1^2 y_3 + q_1 q_2 y_3 \\ 
 &- & q_2 x_1^2 y_3 + 2 q_1 x_1 x_2 y_3 + 
  x_1^2 x_2^2 y_3 + q_1 x_1 y_1 y_3 -q_2 x_1 y_1 y_3\\ 
 & +& q_1 x_2 y_1 y_3 + 
  x_1^2 x_2 y_1 y_3 + x_1 x_2^2 y_1 y_3 + 
  q_1 y_1^2 y_3 +x_1 x_2 y_1^2 y_3\\ 
 &+ & q_1 x_1 y_2 y_3 - q_2 x_1 y_2 y_3 + 
  q_1 x_2 y_2 y_3 + x_1^2 x_2 y_2 y_3 +x_1 x_2^2 y_2 y_3\\ 
 &- & q_2 y_1 y_2 y_3 + x_1^2 y_1 y_2 y_3 + 
  2 x_1 x_2 y_1 y_2 y_3 + x_2^2 y_1 y_2 y_3 + x_1 y_1^2 y_2 y_3 \\ 
 &+ & x_2 y_1^2 y_2 y_3 + q_1 y_2^2 y_3 + 
  x_1 x_2 y_2^2 y_3 + x_1 y_1 y_2^2 y_3 + x_2 y_1 y_2^2 y_3 \\ 
 &+ & y_1^2 y_2^2 y_3, \\ 
\wt{\s}_{12321}(x,y) & = & -q_1 q_2 x_1 + q_2 x_1^3 - 
q_1^2 x_3 + q_1 x_1^2 x_3 - q_1 x_1 x_2 x_3 + x_1^3 x_2 x_3\\ 
 &- & q_1^2 y_1 - q_1 q_2 y_1 + 
q_1 x_1^2 y_1 + q_2 x_1^2 y_1 - q_1 x_1 x_2 y_1 + 
   x_1^3 x_2 y_1\\ 
 &- & q_1 x_1 x_3 y_1 + x_1^3 x_3 y_1 - 
q_1 x_2 x_3 y_1 + x_1^2 x_2 x_3 y_1 - 
   q_1 x_1 y_1^2 + x_1^3 y_1^2\\ 
 &- & q_1 x_2 y_1^2 + 
x_1^2 x_2 y_1^2 - q_1 x_3 y_1^2 + 
   x_1^2 x_3 y_1^2 - q_1 y_1^3 + x_1^2 y_1^3 + 
q_2 x_1^2 y_2\\ 
 &+ & q_1 x_1 x_3 y_2 + 
   x_1^2 x_2 x_3 y_2 + q_1 x_1 y_1 y_2 + q_2 x_1 y_1 y_2 + x_1^2 x_2 y_1 y_2\\ 
 &+ & x_1^2 x_3 y_1 y_2 + x_1 x_2 x_3 y_1 y_2 + x_1^2 y_1^2 y_2 + 
x_1 x_2 y_1^2 y_2 +  x_1 x_3 y_1^2 y_2 \\ 
 &+ & x_1 y_1^3 y_2 + q_2 x_1^2 y_3 + 
q_1 x_1 x_3 y_3 + x_1^2x_2 x_3 y_3 + 
   q_1 x_1 y_1 y_3 \\ 
 &+ & q_2 x_1 y_1 y_3 + 
x_1^2 x_2 y_1 y_3 + x_1^2 x_3 y_1 y_3 + 
   x_1 x_2 x_3 y_1 y_3 + x_1^2 y_1^2 y_3 \\ 
 &+ & x_1 x_2 y_1^2 y_3 + x_1 x_3 y_1^2 y_3 + 
 x_1 y_1^3 y_3 + q_2 x_1 y_2 y_3 + q_1 x_3 y_2 y_3 \\ 
 &+ & x_1 x_2 x_3 y_2 y_3 + q_1 y_1 y_2 y_3 + 
 q_2 y_1 y_2 y_3 + x_1 x_2 y_1 y_2 y_3 \\
 &+ & x_1 x_3 y_1 y_2 y_3 + x_2 x_3 y_1 y_2 y_3 + 
 x_1 y_1^2 y_2 y_3 + x_2 y_1^2 y_2 y_3  \\
&+& x_3 y_1^2 y_2 y_3 + y_1^3 y_2 y_3,  \\ 
\wt{\s}_{12132}(x,y) & = & q_1 q_2 x_1 + q_2 x_1^2 x_2 + 
q_1^2 x_3 + 2 q_1 x_1 x_2 x_3 + x_1^2 x_2^2 x_3 + 
q_1^2 y_1\\ 
 &+ &  q_1 q_2 y_1 + 2 q_1 x_1 x_2 y_1 + 
q_2 x_1 x_2 y_1 + x_1^2 x_2^2 y_1 + 
 q_1 x_1 x_3 y_1\\ 
 &+ & q_1 x_2 x_3 y_1 + 
x_1^2 x_2 x_3 y_1 + x_1 x_2^2 x_3 y_1 + q_1 x_1 y_1^2 + 
 q_1 x_2 y_1^2\\ 
 &+ & x_1^2 x_2 y_1^2 + 
x_1 x_2^2 y_1^2 + q_1 x_3 y_1^2 + x_1 x_2 x_3 y_1^2 + 
   q_1 y_1^3 + x_1 x_2 y_1^3\\ 
 &+ & q_2 x_1^2 y_2 + q_2 x_1 x_2 y_2 + q_1 x_1 x_3 y_2 + 
   q_1 x_2 x_3 y_2 + x_1^2 x_2 x_3 y_2\\ 
 &+ & x_1 x_2^2 x_3 y_2 + q_1 x_1 y_1 y_2 + q_2 x_1 y_1 y_2 + 
   q_1 x_2 y_1 y_2 + q_2 x_2 y_1 y_2\\ 
 &+ & x_1^2 x_2 y_1 y_2 + x_1 x_2^2 y_1 y_2 + x_1^2 x_3 y_1 y_2 + 
2 x_1 x_2 x_3 y_1 y_2 + x_2^2 x_3 y_1 y_2\\ 
 &+ & x_1^2 y_1^2 y_2 + 
   2 x_1 x_2 y_1^2 y_2 + x_2^2 y_1^2 y_2 + x_1 x_3 y_1^2 y_2 + 
x_2 x_3 y_1^2 y_2\\ 
 &+ &  x_1 y_1^3 y_2 + x_2 y_1^3 y_2 + q_2 x_1 y_2^2 + q_1 x_3 y_2^2 + 
x_1 x_2 x_3 y_2^2 + q_1 y_1 y_2^2\\ 
 &+ & q_2 y_1 y_2^2 + x_1 x_2 y_1 y_2^2 + x_1 x_3 y_1 y_2^2 + 
x_2 x_3 y_1 y_2^2 + 
   x_1 y_1^2 y_2^2\\ 
 &+ & x_2 y_1^2 y_2^2 + x_3 y_1^2 y_2^2 + y_1^3 y_2^2, \\ 
\wt{\s}_{1321}(x,y) & = & - q_1^2 - q_1 q_2 + q_1 x_1^2 - 
  q_1 x_1 x_2 + x_1^3 x_2 - q_1 x_1 y_1 + x_1^3 y_1 \\ 
 &- &  q_1 x_2 y_1 + x_1^2 x_2 y_1 - 
  q_1 y_1^2 +x_1^2 y_1^2 + 
  q_1 x_1 y_2 + x_1^2 x_2 y_2 \\ 
 &+ &  x_1^2 y_1 y_2 + 
  x_1 x_2 y_1 y_2 + x_1 y_1^2 y_2 + q_1 x_1 y_3 + 
  x_1^2 x_2 y_3 + x_1^2 y_1 y_3\\ 
 &+ &  x_1 x_2 y_1 y_3 + x_1 y_1^2 y_3 + q_1 y_2 y_3 + 
  x_1 x_2 y_2 y_3 + x_1 y_1 y_2 y_3\\ 
 &+ & x_2 y_1 y_2 y_3 + y_1^2 y_2 y_3, \\ 
\wt{\s}_{2321}(x,y) & = & - q_1^2 - q_1 q_2 + q_1 x_1^2 - 
  q_1 x_1 x_2 + x_1^3 x_2 - 2 q_1 x_1 x_3 + x_1^3 x_3 \\ 
 &- & q_1 x_2 x_3 - q_1 x_1 y_1 + 
  x_1^3 y_1 - q_1 x_2 y_1 + 
  x_1^2 x_2 y_1 - q_1 x_3 y_1\\ 
 &+ &   x_1^2 x_3 y_1 - q_1 y_1^2 + x_1^2 y_1^2 - q_1 x_1 y_2 + 
  x_1^3 y_2 - q_1 x_2 y_2 +x_1^2 x_2 y_2\\  
 &- & q_1 x_3 y_2 + 
  x_1^2 x_3 y_2 - q_1 y_1 y_2 + 
  2 x_1^2 y_1 y_2 + 
  x_1 x_2 y_1 y_2 + x_1 x_3 y_1 y_2\\ 
 &+ & x_1 y_1^2 y_2 - q_1 y_2^2 + 
  x_1^2 y_2^2 + x_1 y_1 y_2^2 + q_1 x_1 y_3 + x_1^2 x_2 y_3 - 
  q_1 x_3 y_3\\ 
 &+ & x_1^2 x_3 y_3 + 
  x_1^2 y_1 y_3 + x_1 x_2 y_1 y_3 + 
  x_1 x_3 y_1 y_3 + 
  x_1 y_1^2 y_3\\ 
 &+ & x_1^2 y_2 y_3 + 
  x_1 x_2 y_2 y_3 + x_1 x_3 y_2 y_3 + 
  2 x_1 y_1 y_2 y_3 + 
  x_2 y_1 y_2 y_3\\ 
 &+ & x_3 y_1 y_2 y_3 + 
  y_1^2 y_2 y_3 + x_1 y_2^2 y_3 + 
  y_1 y_2^2 y_3, \\ 
\wt{\s}_{2132}(x,y) & = & q_1^2 + q_1 q_2 - q_2 x_1^2 + 
  2 q_1 x_1 x_2 + x_1^2 x_2^2 + q_1 x_1 y_1 - q_2 x_1 y_1\\ 
 &+ & q_1 x_2 y_1 + x_1^2 x_2 y_1 + 
  x_1 x_2^2 y_1 + q_1 y_1^2 + x_1 x_2 y_1^2 + q_1 x_1 y_2 \\ 
 &- & q_2 x_1 y_2 + q_1 x_2 y_2 + 
  x_1^2 x_2 y_2 + x_1 x_2^2 y_2 - 
  q_2 y_1 y_2 + x_1^2 y_1 y_2\\ 
 &+ & 2 x_1 x_2 y_1 y_2 + 
  x_2^2 y_1 y_2 + x_1 y_1^2 y_2 + 
  x_2 y_1^2 y_2 + q_1 y_2^2 +x_1 x_2 y_2^2 \\ 
 &+ & x_1 y_1 y_2^2 + x_2 y_1 y_2^2 + y_1^2 y_2^2, \\ 
\wt{\s}_{1213}(x,y) & = & q_2 x_1^2 + q_1 x_1 x_3 + 
  x_1^2 x_2 x_3 + q_1 x_1 y_1 + 
  q_2 x_1 y_1 + x_1^2 x_2 y_1\\ 
 &+ & x_1^2 x_3 y_1 + 
  x_1 x_2 x_3 y_1 + x_1^2 y_1^2 + 
  x_1 x_2 y_1^2 + x_1 x_3 y_1^2 + x_1 y_1^3 \\ 
 &+ & q_2 x_1 y_2 + 
  q_1 x_3 y_2 + x_1 x_2 x_3 y_2 +  q_1 y_1 y_2 + q_2 y_1 y_2 + 
  x_1 x_2 y_1 y_2\\ 
 &+ & x_1 x_3 y_1 y_2 + 
  x_2 x_3 y_1 y_2 + x_1 y_1^2 y_2 + x_2 y_1^2 y_2 + 
  x_3 y_1^2 y_2 + y_1^3 y_2,\\ 
\wt{\s}_{1232}(x,y) & = & q_2 x_1^2 + q_2 x_1 x_2 + 
  q_1 x_1 x_3 + q_1 x_2 x_3 + 
  x_1^2 x_2 x_3 + x_1 x_2^2 x_3\\ 
 &+ & q_1 x_1 y_1 + q_2 x_1 y_1 + 
  q_1 x_2 y_1 + q_2 x_2 y_1 + x_1^2 x_2 y_1 + x_1 x_2^2 y_1\\ 
 &+ &  x_1^2 x_3 y_1 + 2 x_1 x_2 x_3 y_1 + 
  x_2^2 x_3 y_1 + x_1^2 y_1^2 + 
  2 x_1 x_2 y_1^2 + x_2^2 y_1^2\\ 
 &+ &  x_1 x_3 y_1^2 + x_2 x_3 y_1^2 + x_1 y_1^3 + x_2 y_1^3 + 
  q_2 x_1 y_2 + q_1 x_3 y_2  \\ 
 &+ & x_1 x_2 x_3 y_2 + q_1 y_1 y_2 + 
  q_2 y_1 y_2 + x_1 x_2 y_1 y_2 + 
  x_1 x_3 y_1 y_2  \\ 
 &+ & x_2 x_3 y_1 y_2 + x_1 y_1^2 y_2 + x_2 y_1^2 y_2 + 
  x_3 y_1^2 y_2 + y_1^3 y_2 + 
  q_2 x_1 y_3  \\ 
 &+ & q_1 x_3 y_3 + 
  x_1 x_2 x_3 y_3 + q_1 y_1 y_3 + 
  q_2 y_1 y_3 + x_1 x_2 y_1 y_3 + \\ 
 &+ &  x_1 x_3 y_1 y_3 + 
  x_2 x_3 y_1 y_3 + 
  x_1 y_1^2 y_3 + x_2 y_1^2 y_3 + 
  x_3 y_1^2 y_3 + y_1^3 y_3, \\ 
\wt{\s}_{121}(x,y) & = & q_1 x_1 + x_1^2 x_2 + x_1^2 y_1 + 
  x_1 x_2 y_1 + x_1 y_1^2 + 
  q_1 y_2 + x_1 x_2 y_2 \\ 
 &+ & x_1 y_1 y_2 + x_2 y_1 y_2 + 
  y_1^2 y_2, \\
\wt{\s}_{132}(x,y) & = & q_1 x_1 - q_2 x_1 + q_1 x_2 + 
  x_1^2 x_2 + x_1 x_2^2 - q_2 y_1 + x_1^2 y_1 + 2 x_1 x_2 y_1\\ 
 &+ &  x_2^2 y_1 + x_1 y_1^2 + 
  x_2 y_1^2 + q_1 y_2 +  x_1 x_2 y_2 + x_1 y_1 y_2 + 
  x_2 y_1 y_2\\ 
 &+ & y_1^2 y_2 + 
  q_1 y_3 +  x_1 x_2 y_3 + x_1 y_1 y_3 + x_2 y_1 y_3 + 
  y_1^2 y_3, \\
\wt{\s}_{232}(x,y) & = & q_1 x_1 + q_1 x_2 + q_2 x_2 + 
 x_1^2 x_2 + x_1 x_2^2 - q_1 x_3 +  x_1^2 x_3 + x_1 x_2 x_3\\ 
 &+ & x_2^2 x_3 + x_1^2 y_1 + 
 2 x_1 x_2 y_1 + x_2^2 y_1 + x_1 x_3 y_1 + x_2 x_3 y_1 + 
 x_1 y_1^2\\ 
 &+ & x_2 y_1^2 + 
 x_1^2 y_2 + 2 x_1 x_2 y_2 + x_2^2 y_2 + x_1 x_3 y_2 + 
  x_2 x_3 y_2 + 2 x_1 y_1 y_2\\ 
 &+ &  2 x_2 y_1 y_2 + x_3 y_1 y_2 +y_1^2 y_2 + x_1 y_2^2 + 
  x_2 y_2^2 + y_1 y_2^2 + q_1 y_3 + q_2 y_3 \\ 
 &+ &  x_1 x_2 y_3 + 
  x_1 x_3 y_3 + x_2 x_3 y_3 + 
  x_1 y_1 y_3 + x_2 y_1 y_3 + 
  x_3 y_1 y_3 \\ 
 &+ & y_1^2 y_3 + 
  x_1 y_2 y_3 + x_2 y_2 y_3 + 
  x_3 y_2 y_3 + y_1 y_2 y_3 + 
  y_2^2 y_3, \\ 
\wt{\s}_{123}(x,y) & = & q_2 x_1 +  q_1 x_3 + x_1 x_2 x_3 + 
  q_1 y_1 + q_2 y_1 +  x_1 x_2 y_1 +x_1 x_3 y_1 \\ 
 &+ & x_2 x_3 y_1 + 
   x_1 y_1^2 + x_2 y_1^2 + 
  x_3 y_1^2 + y_1^3, \\ 
\wt{\s}_{213}(x,y) & = & q_1 x_1 + x_1^2 x_2 - q_1 x_3 + 
  x_1^2 x_3 + x_1^2 y_1 + 
  x_1 x_2 y_1 +x_1 x_3 y_1 + 
  x_1 y_1^2 \\ 
 &+ & x_1^2 y_2 + 
  x_1 x_2 y_2 + x_1 x_3 y_2 + 
  2 x_1 y_1 y_2 +x_2 y_1 y_2 + 
   x_3 y_1 y_2 + y_1^2 y_2 \\ 
 &+ & x_1 y_2^2 +  y_1 y_2^2, \\ 
\wt{\s}_{321}(x,y) & = & -2 q_1 x_1 + x_1^3 - q_1 x_2 - 
  q_1 y_1 + x_1^2 y_1 - q_1 y_2 +x_1^2 y_2 + x_1 y_1 y_2 \\ 
 &- &  q_1 y_3 + x_1^2 y_3 + 
  x_1 y_1 y_3 + x_1 y_2 y_3 + y_1 y_2 y_3, \\ 
\wt{\s}_{23}(x,y) & = & q_1 + q_2 + x_1 x_2 + x_1 x_3 + 
  x_2 x_3 +  x_1 y_1 + x_2 y_1 +x_3 y_1 + y_1^2 \\ 
 &+ & x_1 y_2 + 
  x_2 y_2 + x_3 y_2 + y_1 y_2 + 
  y_2^2, \\ 
\wt{\s}_{32}(x,y) & = & - q_1 - q_2 + x_1^2 + 
x_1 x_2 + x_2^2 + 
x_1 y_1 + x_2 y_1 + x_1 y_2 + 
  x_2 y_2 \\ 
 &+ & y_1 y_2 + x_1 y_3 + 
  x_2 y_3 + y_1 y_3 + y_2 y_3, \\ 
\wt{\s}_{13}(x,y) & = &  x_1^2 + x_1 x_2 + x_1 x_3 + 
  2 x_1 y_1 + x_2 y_1 + x_3 y_1 + y_1^2 + x_1 y_2  + y_1 y_2\\ 
 &+ &   x_1 y_3 + y_1 y_3, \\ 
\wt{\s}_{12}(x,y) & = &  q_1 + x_1 x_2 + x_1 y_1 + x_2 y_1 + 
  y_1^2, \\ 
\wt{\s}_{21}(x,y) & = & - q_1 + x_1^2 + x_1 y_1 + x_1 y_2 + 
  y_1 y_2, \\ 
\wt{\s}_{3}(x,y) & = & x_1 + x_2 + x_3 + y_1 + y_2 + y_3, \\ 
\wt{\s}_{2}(x,y) & = & x_1 + x_2 + y_1 + y_2, \\ 
\wt{\s}_{1}(x,y) & = & x_1 + y_1, \\ 
\wt{\s}_{id}(x,y) & = & 1.
\end{eqnarray*} 

\section*{Appendix B}
{\large{\bf Quantum residue generating function $\Psi (t)$ for $S_3$.}}

\vskip 0.5cm
We define the functions 
$g_{v}(t_1)$ and $h_{v}(t_2)$ by the formulas 
\[ \frac{t_1}{t_1-x_1} = \sum_{v\in S_3}
g_{v}(t_1){\cal S}^q_v, \]
\[ \frac{t_2}{t_2-x_2} = \sum_{v\in S_3} 
h_{v}(t_2){\cal S}^q_v \] 
in the quantum cohomology ring $QH^{\ast}(Fl_3).$
Then, we have 
\[ \langle \frac{t_1}{t_1 - x_1} 
\frac{t_2}{t_2 - x_2} \rangle 
= \sum_{v\in S_3} g_v (t_1) h_{vw_0} (t_2). \] 
Since $\langle x_1^5 \rangle = q_1,$ 
$x_2 = q_1^{-1} x_1^3 - 2x_1$ and 
\[ QH^{\ast}(Fl_3) \simeq {\bf Z}[q_1, q_2][x_1]/(f_1(x_1)), \] 
the correlation function $\langle P(x_1,x_2) \rangle $ is 
expressed as 
$$\langle P(x_1,x_2) \rangle  =  
q_1 {\rm Res}_{f_1} P(x_1,q_1^{-1}x_1^3-2x_1) 
= q_1 \sum_{f_1(\mu)=0}\frac{1}{f'_1(\mu)}
P(\mu,q_1^{-1}\mu^3-2\mu) .
$$
Similarly, it also holds that 
\[ \langle P(x_1,x_2) \rangle = 
(q_2-q_2) \sum_{f_2(\mu)=0} 
\frac{1}{f'_2(\mu)} P((q_1-q_2)^{-1}(\mu^3-(2q_1+ q_2)\mu),
\mu). \] 
The functions $g_v$ and $h_v$ are given as follows: 
\begin{eqnarray*} 
f_1(t_1) g_{121}(t_1) & = & q_1t_1, \\ 
f_1(t_1) g_{12}(t_1) & = & q_1t_1^2, \\  
f_1(t_1) g_{21}(t_1) & = & t_1^2(t_1^2-q_1), \\ 
f_1(t_1) g_{2}(t_1) & = & q_1t_1(t_1^2-q_1), \\ 
f_1(t_1) g_1(t_1) & = & t_1(t_1^4-q_1t_1^2+q_1^2), \\ 
f_1(t_1) g_{id}(t_1) & = &  t_1^2(t_1^4-q_1t_1^2+q_1^2), \\ 
f_2(t_2) h_{121}(t_2) & = & (q_2-q_1)t_2, \\ 
f_2(t_2) h_{12}(t_2) & = & t_2^2(2q_1+q_2-t_2^2), \\
f_2(t_2) h_{21}(t_2) & = & t_2^2(q_1+2q_2-t_2^2), \\ 
f_2(t_2) h_2(t_2) & = &  \frac{(q_1+q_2)t_2^5-q_1(q_1+q_2)t_2^3+q_1q_2(2q_2-q_1)t_2}
{(q_1-q_2)}, \\ 
f_2(t_2) h_1(t_2) & = & 
-t_2^5+(2q_1+q_2)t_2^3+q_1(q_2-q_1)t_2, \\ 
f_2(t_2) h_{id}(t_2) & = & \frac{(q_1+q_2)t_2^6-2q_1t_2^2+
(q_1^3+2q_1^2q_2-q_2^3)t_2^2}{(q_1-q_2)}. 
\end{eqnarray*} 
The characteristic polynomials $f_1$ and $f_2$ are given by 
\[ f_1 (t) = ( t^2-q_1)^3 -q_1^2 q_2, \] 
\[ f_2 (t) = t^6 - 3 (q_1 + q_2)t^4 + 3(q_1^2 + q_1 q_2 + q_2^2) 
t^2 -q_1^3 + q_1^2q_2 + q_1 q_2^2 -q_2^3 . \] 
Hence, we have 
\begin{eqnarray*} 
&D_1& =  \frac{1}{720t_1}f_1(t_1)
\frac{\partial^6}{\partial t_1^6} + 
\left( \frac{q_1^3 + q_1^2 q_2}{120t_1^2} 
+\frac{1}{40}q_1^2 -\frac{3}{40}q_1t_1^2 
+\frac{1}{24}t_1^4 \right) 
\frac{\partial^5}{\partial t_1^5}  \\ 
&+& \left( -\frac{q_1^3 + q_1^2 q_2}{24t_1^3} 
- \frac{3}{8}q_1 t_1 + \frac{5}{12}t_1^3 \right) 
\frac{\partial^4}{\partial t_1^4} + 
\left( \frac{q_1^3 + q_1^2 q_2}{6t_1^4} 
-\frac{1}{2}q_1 +\frac{5}{3}t_1^2 \right) 
\frac{\partial^3}{\partial t_1^3}  \\ 
&+& \left( -\frac{q_1^3 + q_1^2 q_2}{2t_1^5} + 
\frac{5}{2}t_1 \right) 
\frac{\partial^2}{\partial t_1^2} + 
\left( \frac{q_1^3 + q_1^2 q_2}{t_1^6} + 1 \right) 
\frac{\partial}{\partial t_1} 
-\frac{q_1^3 + q_1^2 q_2}{t_1^7}, \\ 
&D_2&=\frac{f_2(t_2)}{720t_2} 
\frac{\partial^6}{\partial t_2^6}  \\ 
&+& \left( \frac{q_1^3 - q_1^2 q_2 - q_1 q_2^2 + q_2^3}{120 t_2^2} 
+ \frac{q_1^2 + q_1q_2 + q_2^2}{40} + 
 \frac{3}{40}(q_1 + q_2) t_2^2 + \frac{t_2^4}{24} \right) 
\frac{\partial^5}{\partial t_2 ^5}  \\ 
&+& \left( -\frac{q_1^3 - q_1^2 q_2 - q_1 q_2^2 + q_2^3}
{24 t_2^3} 
- \frac{3}{8}(q_1 + q_2) t_2 + \frac{5}{12}t_2^3 \right) 
\frac{\partial^4}{\partial t_2 ^4}  \\ 
&+& \left( \frac{q_1^3 - q_1^2 q_2 - q_1 q_2^2 + q_2^3}{6 t_2^4} 
-\frac{1}{2}(q_1 + q_2) + \frac{5}{3}t_2^2 \right) 
\frac{\partial^3}{\partial t_2 ^3}  \\ 
&+& \left( -\frac{q_1^3 - q_1^2 q_2 - q_1 q_2^2 + q_2^3}{2t_2^5} 
+ \frac{5}{2}t_2 \right) 
\frac{\partial^2}{\partial t_2 ^2}  \\ 
&+& \left( \frac{q_1^3 - q_1^2 q_2 - q_1 q_2^2 + q_2^3}{t_2^6} 
+ 1 \right) 
\frac{\partial}{\partial t_2}  
 - \frac{q_1^3 - q_1^2 q_2 - q_1 q_2^2 + q_2^3}{t_2^7} . 
\end{eqnarray*} 



\vskip 1cm

\begin{thebibliography}{999}
        \bibitem [BB]{[BB]}
        Bergeron N. and Billey S., {\it RC--graphs and Schubert polynomials,}
        Experimental Math. 2, 1993, p.257-269;

        \bibitem [BJS]{[BJS]}
        Billey S., Jockusch W., Stanley R., {\it Some combinatorial properties 
        of Schubert polynomials}, Journal of Algebraic Combinatorics, 1993, 
        v.2, 345-374; 
        
	\bibitem [C]{[C]}
	Ciocan--Fontanine I., {\it Quantum cohomology of flag varieties},
        Intern. Math.Research Notes, 1995, n.6, p.263-277;
        
	\bibitem [EL]{[EL]}
	Eisenbud D., Levine H., {\it An algebraic formula for the degree
        of a $C^{\infty}$ map germ}, Ann. of Math., 1977, v.106, p.19-38;
        
	\bibitem [Fl]{[FI]}
	Di Francesco P., Itzykson C., {\it Quantum intersection rings, 
        in} The Moduli Space of Curves, Progress in Mathematics, 
        Birkh\"auser, 1995, v.129, p.81-148; 
        
        \bibitem [F]{[F]}
        Fulton W., {\it Flags, Schubert polynomials, degeneracy loci, and 
        determinantal formulas}, Duke Math. J., 1991, v.65, p.381-420;

        \bibitem [FGP]{[FGP]}
        Fomin S, Gelfand S. and Postnikov A., {\it Quantum Schubert
        polynomials,} Preprint, 1996, 44p.;
        
	\bibitem [GK]{[GK]}	
        Givental A., Kim B., {\it Quantum cohomology of flag manifolds and 
        Toda lattices}, Comm.Math.Phys., 1995, v.168, p.609-641;
        
        \bibitem [GH]{[GH]}
        Griffiths Ph., Harris J., {\it Principles of algebraic geometry}, 
        John Wiley\& Sons, New York, 1978.
        
	\bibitem [K1]{[K1]}
	Kim B., {\it Quot schemes for flags and Gromov invariants for flag 
	varieties}, Prepr., 1995, alg-geom/9512003;
	
	\bibitem [K2]{[K2]}
	Kim B., {\it On equivariant quantum cohomology}, Prepr., 1995, 
	\hbox{q-alg/9509029;}
	
	\bibitem [K3]{[K3]}
	Kim B., {\it Quantum cohomology of flag manifolds $G/B$ and quantum 
        Toda lattices}, 1996, Prepr., alg-geom/9607001;
        
        \bibitem [KV]{[KV]}
        Kohnert A., Veigneau S., {\it Using Schubert toolkit to compute with 
        polynomials in several variables}, in 8-th International Conf. on 
        Formal Power Series and Algebraic Combinatorics, Univ. of Minnesota, 
        1996, p.283-293 ;
        
	\bibitem [KM]{[KM]}
	Kontsevich M., Manin Yu., {\it Gromov-Witten classes, quantum 
        cohomology, and enumerative geometry}, Comm.Math.Phys., 1994,
        v.164, p.525-562;

	\bibitem [LS1]{[LS1]}
	Lascoux A., Sch\"utzenberger M.-P., {\it Polyn\^omes de Schubert,}
        C.R. Acad. Sc. Paris, t.294, 1982, p.447-450;
        
	\bibitem [LS2]{[LS2]}
	Lascoux A., Sch\"utzenberger M.-P., {\it Symmetry and flag manifolds}, 
	Lect. Notes in Math., 1983, v.996, p.118-144;
        
	\bibitem [M1]{[M1]}
	Macdonald I.G., {\it Notes on Schubert polynomials}, Publ. LCIM, 1991, 
	Univ. de Quebec a Montreal;

       \bibitem [M2]{[M2]}
        Macdonald I.G., {\it Symmetric Functions and Hall Polynomials},
        second edition, Oxford University Press, new York/London, 1995;
	
	\bibitem [MS]{[MS]}
	McDuff D., Salamon D., {\it $J$-holomorphic curves and quantum 
	cohomology}, Univ. Lect., v.6, 1994, AMS;

       \bibitem [RT]{[RT]}
        Ruan Y., Tian G., {\it Mathematical theory of quantum cohomology},
        J.Diff.Geom., 1995, v.42, n.2, p.259-367;

        \bibitem [S]{[S]}
        Sotile F., {\it Pieri's rule for flag manifolds and Schubert
        polynomials,} Preprint, 1995, 20p.


\end{thebibliography}
\end{document}